\documentclass{article}

\usepackage{amsmath, amsfonts,amssymb, amsbsy}

\newtheorem{theo}{\bf Theorem}[section]
\newtheorem{cor}{\bf Corollary}[section]
\newtheorem{defi}{\bf Definition}[section]
\newtheorem{nota}{\bf Notation}[section]
\newtheorem{lem}{\bf Lemma}[section]

\newcommand{\ba}{\boldmath$\alpha$}

\begin{document}

\begin{center}
{\bf {\Large Inter-class orthogonal main effect plans for\\
asymmetrical experiments.}}
\end{center}
\begin{center}
  Sunanda Bagchi,\\
   Theoretical Statistics and Mathematics Unit,
Indian Statistical Institute,\\ Bangalore 560059, India.
\end{center}
\vskip5pt {\bf {\large Abstract : }} In this paper we construct
`inter-class orthogonal' main effect plans (MEP) for asymmetrical
experiments. In such a plan, a factor is orthogonal to all others
except possibly the ones in its own class. We have also defined the
concept of ``partial orthogonality" between  a pair of factors. In
many of our plans, ``partial orthogonality"  has been achieved when
(total) orthogonality is not possible  due to divisibility or any
other restriction.

We present a method of obtaining `inter-class orthogonal' MEPs.
Using this method and also a method of `cut and paste' we have
obtained several series of `inter-class orthogonal' MEPs.
Interestingly some of these happen to be  orthogonal MEP (OMEP); for
example we have constructed  an OMEP for a $3^{30}$ experiment on 64
runs. Further, many of the `inter-class orthogonal' MEPs are `almost
orthogonal' in the sense that each factor is orthogonal to all
others except possibly one. In many of the other MEPs factors are
``orthogonal through another factor", thus leading to simplification
in the analysis. Plans of small size ($\leq 15$ runs) are also
constructed by ad-hoc methods.

Finally, we present a user-friendly computational method for
analysing data obtained from any general factorial design.



\section{Introduction}

 In many industrial experiments like screening experiments,
 often the interest lies only in the main effects of
factors. The wide use of orthogonal main effect plans (OMEP) for
such experiments is due to their orthogonality property which
ensures uncorrelated and hence most precise estimation of every main
effect contrast of every factor, apart from providing great
simplicity in analysis, as is well-known.

However, orthogonality requires certain divisibility conditions and
so an OMEP for an asymmetrical experiment often requires a large
number of runs. [See Dey  and Mukherjee  (1999) and Hedayat, Sloan,
and Stuffken (1999) for details]. The proportional
  frequency (PF) plans of Addleman (1962), as we know, are  OMEPs
 with possibly unequal replications for one or more factors, thus
 requiring weaker conditions for existence. However, very few
 unequally replicated PF plans are known, apart from
 Stark's (1964) plan for a $3^7$ experiment on 16 runs.
 Thus, the problem of availability of an OMEP with not-too-large
 run size remains. In such situations, therefore, the question
 arises whether with a smaller run size one can find an alternative
 plan  - something not as good as an OMEP but not too bad either.

     Of late, departure from full orthogonality has been investigated
in the context of main effect plans (MEPs). In the  ``nearly
orthogonal" plans of Wang and Wu (1992) factors are allowed to be
non-orthogonal to a few of the other factors. Subsequently, other
nearly orthogonal MEPs having interesting  combinatorial properties
have been proposed and studied by others like Nguyen (1996), Ma,
Fang and Liski (2000), Huang, Wu and Yen (1992) and  Xu (2002).

       Why do we look for ``near orthogonality"? Why can't we go far
       away and use a fully non-orthogonal plan ? If we are willing
       to use non-orthogonal plans, we would have tremendous
       flexibility. We could, for instance, make a $2^4$ experiment on 5
       runs (instead of 8)  [see plan $A_5(4)$ of example 2.1] or a $3^5$ experiment on
       12 runs (instead of 16) [see Plan $A_{12}(4)$ in section 5]. One hurdle
       to the usability of such plans is the complexity in the
 data analysis. The reduction in the precision is, of course, another problem.

  In the present paper our main aim is to provide main effect
plans (MEPs) for asymmetrical experiments with small run size,
deviating ``as little as possible" from the desirable properties
like orthogonality and/or equal replications, so that analysis
remains relatively simple. Specifically, we construct plans
satisfying {\bf ``inter-class orthogonality"}, in which each factor
is possibly non-orthogonal to the members of its own class, but
orthogonal to factors of other classes. In the process we have also
obtained a series of {\bf orthogonal MEP  for a $3^{30}$ experiment
on 64 runs} (see Theorem \ref {two-stage-Stark}). In many of our
plans the class size is at most two, so that a factor is orthogonal
to all others except possibly one. Among plans of larger class size,
many plans satisfy the property that within class factors are
``orthogonal through another factor" (in the same class), thus
leading to simplification in the analysis. (See Example 2.1 and
Theorem \ref {two-stage} (c)).

We have also defined the concept of {\bf ``partial orthogonality"}
between a pair of factors and derived sufficient condition for it.
[See definition \ref {partOrth}, Lemma \ref {condPFC} and the
discussion thereafter.  In many of our plans, ``partial
orthogonality" has been achieved between one or more pair(s) of
factors when (full) orthogonality is
 not possible  due to divisibility or any other restriction.

   The definitions along with examples are presented in section 2.
   In section 3 we construct a few series of  ``inter-class orthogonal" MEPs
for asymmetrical experiments with factors having at most five
levels.  Using ad-hoc methods we have also constructed MEPs with
factors nonorthogonal to one or two factors on at most 15 runs,
which are in Section 5. These plans include saturated plans for the
following experiments. $4^2.2$, $3^2.2^3$ and $4.3.2^2$  on 8 runs,
$5^2.2$ on 10 runs, $4^2.3^2.2$ and $2.3^5$ on 12 runs and $5^2.3^2$
on 15 runs. In section 4 we
   present an user-friendly method of analysis.

 We believe that the information presented in Theorem \ref {InclAna}
 and other results in section 4 will help the experimenter to have
  a clear idea about the efficiencies of the
BLUEs of the main effects as well as the amount of computation
involved in the analysis of a non-orthogonal plan. These features
may be compared with those of other available plans like ``plan
orthogonal through one factor" or an ``inter-class orthogonal plan".

    This paper has been posted in aexiv.org/abs/1512.06588.

\section{Inter-class orthogonal main effect plans}
 Throughout this paper we shall be concerned with main effect plans,
 that is plans aiming at gaining information only about the main
 effects, assuming interactions to be negligible.
 In all plans presented henceforth in this paper, {\bf rows represent factors,
   while columns represent runs}.

 Let us consider a main effect plan (MEP) for an experiment with $m$
 factors $A, B, \cdots $ on $n$ runs. Suppose the factor $A$ have $a$ levels,
 $B$ have $b$ levels and so on. Then, the plan will be referred to as an
 m-factor MEP and will be represented by an $m \times n$ array {\boldmath
 $\rho(n,m;a,b, \cdots)$}.

  {\boldmath $r_A (i) $} {\bf will denote
 the number of runs in which factor A is at level i}, while the
 vector {\boldmath $r_A = (r_A (1), \cdots r_A (a))$} will be referred
 to as the {\bf replication vector } of the factor A. For two
 factors A,B the {\bf incidence matrix}  {\boldmath $N_{AB}$} is the
  {\boldmath $a \times b$} matrix with the $(i,j)$th entry
 {\boldmath $n_{AB} (i,j)$} = {\bf the number of runs in which A is
 at level i and B is at level j}. Clearly,  {\boldmath $N_{AA}$}
 {\bf is a diagonal matrix whose diagonal entries are those of}
  {\boldmath $r_A$} in the same order and will sometimes
 be denoted by  {\boldmath $R_A$}.

\begin{defi}
Let us consider an m-factor MEP $\rho$ on n runs. Suppose the set of
factors of $\rho$ can be divided into several classes in such a way
that {\bf every factor is orthogonal to every other from a different
class}. Then $\rho$ is called {\bf inter-class orthogonal}. An
inter-class orthogonal  MEP with $m$ factors divided into $p$
classes  and  the factors in the i-th class having levels  $ s_{i1},
s_{i2}, \cdots$ on $n$ runs will be denoted by  ${\boldmath
\rho(n,m; \{s_{11}. s_{12} \cdots \}. \{ s_{21}. s_{22} \cdots \}
\cdots )}$. A plan with at most $m$ factors in a class may be
referred to as an {\bf inter-class(m) orthogonal} MEP.
\end{defi}

{\bf Remark 2.1:} Any main-effect plan may be looked upon as an
inter-class(m) orthogonal MEP, for some m. For instance, an OMEP may
be viewed as an inter-class (1) orthogonal MEP, while an MEP with p
factors of which no one orthogonal to any other is inter-class(p)
orthogonal. The plan $L'_{18} (3^4 . 2^8)$ of Wang and Wu(1992) is
an inter-class (8) orthogonal MEP, according to the present
terminology, as the 8 two-level factors are mutually non-orthogonal.
We see that the term inter-class(m) orthgonal does not always
display the exact picture as there may be classes with size much
smaller than m, as in the case of $L'_{18} (3^4 . 2^8)$. This term
is informative when the class sizes are close to one another, which
is the case for  the plans constructed here.

{\bf Examples :} We now present two inter-class orthogonal plans
 and along with their graphical representation. Here adjacency represents
 orthogonality. The interpretation of the dotted line between factors  is explained
in Remark 2.4.

{\bf Example 2.1:} $A_8(1) = \rho(8,5; \{3^2\}.\{2^2\}.2) $.
\begin{eqnarray*}  A_8(1) =
 \left [ \begin{array}{ccccccccc}
         A&  0 & 1 & 0 & 2 & 0 & 1 & 0 & 2 \\
         B& 0 & 0 & 1 & 2 & 2 & 1 & 0 & 0 \\
         C& 0 & 0 & 0 & 0 & 1 & 1 & 1 & 1 \\
         D& 0 & 1 & 1 & 1 & 0 & 0 & 1 & 0 \\
         E& 0 & 1 & 1 & 0 & 1 & 0 & 0 & 1 \\
           \end{array} \right ] & \; \;\qquad
           \setlength{\unitlength}{3mm}
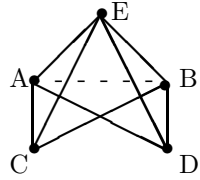
\begin{picture}(20,6)(0,3) \thicklines
\put(1,1){\line(0,1){3}} \put(1,1){\line(2,1){6}}
 \put(1,1){\line(1,2){3}}
 \put(1,4){\line(2,-1){6}} \put(1,4){\line(1,1){3}}
  \put(4,7){\line(1,-2){3}}  \put(4,7){\line(1,-1){3}}
 \put(7,1){\line(0,1){3}} \put(7,1){\line(-1,2){3}}

 \put(0.8,0.8){$\bullet$} \put(0.8,3.8){$\bullet$} \put(3.8,6.8){$\bullet$}
 \put(6.7,0.8) {$\bullet$} \put(6.6,3.6){$\bullet$}

\put(0.0,0.0){C} \put(0.0,3.8){A} \put(4.5,6.8){E} \put(7.5,0.0){D}
\put(7.5,3.8){B} \put(1,3.8){-} \put(2,3.8){-}
\put(3,3.8){-}\put(4,3.8){-}
 \put(5,3.8){-}\put(6,3.8){-}
\end{picture}  \end{eqnarray*}

{\bf Example 2:} $A_{12}(1) =  \rho(12,5; \{2.4^2 \}.\{3^2\})$.
\begin{eqnarray*} A_{12}(1) = \left [ \begin{array}{ccccccccccccc}
       A& 0 & 1 & 0 & 1 & 0 & 1 & 1 & 0 & 1 & 0 & 1 & 0 \\
       B& 0 & 0 & 0 & 1 & 1 & 1 & 2 & 2 & 2 & 3 & 3 & 3 \\
       C& 1 & 2 & 3 & 2 & 3 & 0 & 3 & 0 & 1 & 0 & 1 & 2 \\
       D& 0 & 1 & 2 & 0 & 1 & 2 & 0 & 1 & 2 & 0 & 1 & 2   \\
       E& 0 & 1 & 2 & 2 & 1 & 0 & 0 & 2 & 1 & 1 & 2 & 0 \\
                            \end{array} \right ] & \; \; \qquad
 \setlength{\unitlength}{3mm}
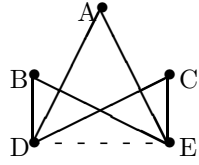
\begin{picture}(20,4)(0,3)

\thicklines

 \put(1,1){\line(0,1){3}} \put(1,1){\line(2,1){6}}
 \put(1,1){\line(1,2){3}}

 \put(1,4){\line(2,-1){6}} \put(4,7){\line(1,-2){3}}
 \put(7,1){\line(0,1){3}}

 \put(0.8,0.8){$\bullet$} \put(0.8,3.8){$\bullet$} \put(3.8,6.8){$\bullet$}
 \put(6.8,0.8) {$\bullet$} \put(6.8,3.8){$\bullet$}

 \put(0.0,0.5){D}  \put(0.0,3.5){B}  \put(3.0,6.5){A}
 \put(7.5,0.5){E}  \put(7.5,3.5){C}

\put(1,0.8){-} \put(2,0.8){-}\put(3,0.8){-}\put(4,0.8){-}
\put(5,0.8){-} \put(6,0.8){-}

 \end{picture} \end{eqnarray*}
 \vspace{.5em}

In this equal frequency saturated plan, both the four-level factors
B and C  form a generalized group divisible design with the levels
of the two-level factor A. Between themselves, they form a Balanced
incomplete block design (BIBD). The relation between factors $D$ and
$E$ is presented in details after Remark 2.4.

\vspace{.5em}

The relation between the factors $A$ and $B$ in $A_8(1)$ and factor
$D$ and $E$ in $A_{12}(1)$ motivates us to define the concept of
partial orthogonality between two factors.

\begin{defi} \label{partOrth}
We say that the factor $A$ is partially orthogonal (PO) to another
factor $B$ if the BLUE of at least one (but not all) main effect
contrast of $A$ is orthogonal to the BLUE of every one of $B$.
 \end{defi}

 {\bf Remark 2.2.} One can verify that in the MEP presented in (2)
of Huang, Wu and Yen (2002), the three-level factors are partially
orthogonal to each other. In fact the relation between every pair of
three-level factors in that plan just like the relation between $A$
and $B$ of $A_8(1)$. [See Table 2.1 below]. More examples are in
Section 3.

{\bf Remark 2.3.} If $A$ is PO to $B$, then $B$ is either PO (plan
$A_8(1)$) to $A$ or non-orthogonal (plan $A_4(2)$ in section 3) to
$A$. Regarding analysis, however, what matters is whether $A$ and
$B$ are mutually orthogonal or not. Thus, partial orthogonality is a
feature of estimation and has no role to play in testing of
hypothesis.

{\bf Remark 2.4.} If two factors are partial orthogonal to each
other, then in the graphical representation they are joined by doted
lines.

{\bf  A statement like ``$A$ is PO to $B$" immediately raise the
question ``which contrast of $A$ is orthogonal to $B$"? We shall now
see how the incidence matrix $N_{AB}$ helps us to find at least
partial answer to this question.}

\vspace{.5em}
\begin{center}
{\bf  How to check orthogonality of a contrast of A to those of B.}
\end{center}

We recall the proportional frequency  condition
 of Addelman (1962).

\begin{defi} \label{ADDL}[ Addelman (1962)] Consider a main effect
plan $\rho$ on $n$ runs.  Two factors $A$ and $B$ are said to be
orthogonal to each other, if the incidence matrix $N_{AB}$ of A and
B satisfies the proportional frequency  condition (PFC), as stated
below.
\begin{equation} \label {PFC}
n_{A,B} (i,j) = r_A (i). r_B (j) /n, \; i=1,2, \cdots a, j =1,2,
\cdots b.\end{equation} \end{defi}

We now define PFC between one factor and certain levels of another
factor.

\begin{defi} \label{partOrth2} Consider two factors $A$ and $B$,
with a and b levels respectively, of a main effect plan $\rho$ on
$n$ runs.

 (a) If a level $i$ of $A$ satisfies
\begin{equation} \label {PFCi}
n_{A,B} (i,j) = r_A (i). r_B (j) /n, \;  j =1,2, \cdots b,
\end{equation}
 then we say that  the level $i$ of $A$ satisfies PFC with factor $B$.

 (b) If a pair of levels $i$ and $k$ of $A$ satisfies
\begin{equation} \label {ikPFC}
n_{A,B} (i,j) / r_A (i) = n_{A,B} (k,j) / r_A (k), \; j =1,2, \cdots
b, \end{equation} then the pair $\{i,k\}$ of levels of $A$ is said
to satisfy PFC with factor $B$.
\end{defi}

   We use the notation {\ba}$_i$ for {\bf the unknown effect of level $i$ of
the factor $A$}, \: $ 1 \leq i \leq a$ and similar notation for
other factors. Further, $\hat{\alpha_i}  - \hat{\alpha_j}$  will
denote the BLUE of the contrast {\ba}$_i -$ {\ba}$_j$. Similar
notation for other contrasts.

The proof of the following result is by straightforward
verification.

\begin{lem} \label{condPFC} (a) If a level $i$ of $A$ satisfy PFC
with $B$, then the BLUE of the  main effect contrast
 $(a-1)${\ba}$_i -   \sum \limits_{j \neq i}$ {\ba}$_{j}$ is orthogonal
 to the BLUEs of all main effect contrasts of $B$.

(b) If  the pair  of levels $\{i,k\}$ of $A$ satisfies PFC with $B$,
then the BLUE of the main effect contrast {\ba}$_i$ - {\ba}$_k$ is
orthogonal to the BLUEs of all main effect contrasts of $B$.
\end{lem}

 We now illustrate these results with the help of plan $A_8(1)$ and
 another plan presented later.

\vspace{.5em}

 {\bf Plan $A_8(1)$ :} We note that factors $A$ and $C$ satisfies PFC
 (see equation (\ref {PFC})) and hence they are mutually orthogonal.
 Similarly, the pairs $(A,D), (A,E),(B,C),(B,D)$ and $(B,E)$ are also mutually orthogonal.
 Regarding the pair of factors $(A,B)B$, we see that PFC condition is not
 satisfied. However, level 0 of $A$ satisfies PFC with factor $B$, as
 shown in the table below.
Therefore, by (a) of Lemma \ref {condPFC} the contrast $
2\hat{\alpha_0} - \hat{\alpha_1} - \hat{\alpha_2}$ is orthogonal to
both the contrasts of $B$. By the same argument the contrast $
2\hat{\beta_0} - \hat{\beta_1} - \hat{\beta_2}$ is orthogonal to
both the contrasts of $A$.

\begin{center}

{\bf {\large Table 2.1}}

 $\left [\begin{array}{cccccccccccc}
 \hline
     &       && B& \rightarrow & &|&  \\\hline
    & A  \downarrow &|&  0& 1 &2&|&r_A \downarrow &|& & r_A.r'_B/n \\\hline
                  &0&|&  2& 1 &1&|& 4             &|&    2 &1   &1  \\
       N_{AB} &    1&|&  1&1  &0&|& 2             &|&    1 & 1/2&1/2   \\
           &       2&|&  1&0  &1&|& 2             &|&    1 & 1/2&1/2 \\\hline
      r_B &  \rightarrow&|    &4&2   &2&|&n=8  \\ \hline
 \end{array}\right ]$.
  \end{center}

 \vspace{.5em}

 {\bf Plan  $A_{12}(1)$ :}
 We note that the pair of three-level factors $D$ and $E$
do not satisfy PFC condition. However, levels 0 and 2 of $D$ satisfy
PFC with $E$ and so by (b) of Lemma \ref  {condPFC} the contrast $
\hat{\delta_0} - \hat{\delta_2}$ is orthogonal to both the contrasts
of $E$. By the same argument the contrast $ \hat{\epsilon_1} -
\hat{\epsilon_2}$ is orthogonal to both the contrasts of $D$. In the
following table $r$ denotes the constant replication number of $D$.

\begin{center}

{\bf {\large Table 2.2}}

 $\left [\begin{array}{cccccccccccc}
 \hline
     &       && E& \rightarrow & &|&  \\\hline
    & D  \downarrow &|&  0& 1 &2&|&r_D \downarrow &|& &  N_{DE}/r \\\hline
                  &0&|&  2& 1 &1&|& 4             &|&    1/2 &1/4   &1/4  \\
       N_{DE} &    1&|&  0&2  &2&|& 4             &|&    0 & 1/2&1/2   \\
           &       2&|&  2&1  &1&|& 4             &|&    1/2 & 1/4&1/4 \\\hline
      r_E & \rightarrow&|&4&4 &4&|&n=12  \\ \hline
 \end{array}\right ]$.
  \end{center}

{\bf Discussion:} What is the use of partial orthogonality ? This
may be viewed as a ``something is better than nothing" approach. If
it is not possible to make $A$ and $B$ mutually orthogonal, we may
at least make them partially orthogonal, if possible. However, the
issue is more complicated, since,  to achieve one condition, we may
have to sacrifice another. Let us look at the following situations.
Consider two factors $A$ and $B$ with $a$ and $b$ levels
respectively.

 Case 1. $ab$ does not divide $n$, so that there does not exist any plan
in which  $A$ and $B$ are mutually orthogonal, each with equal
frequency. In case a proportional frequency plan exists, then of
course, that is the best option. Suppose such a plan is not known.
If we know a plan, say $\rho_1$, in which $A$ is partially
orthogonal to $B$, then the experimenter would be happy to be able
to estimate at least a few among the main effect contrasts of $A$
with maximum precision. However, this may lead to ``too small" a
precision for the other contrasts of $A$. Suppose a plan $\rho_2$ is
also available in which $A$ and $B$ are not partially orthogonal,
but all the contrasts of $A$ and $B$ are estimated with ``reasonably
high" precision. Whether the experimenter would prefer $\rho_1$ or
$\rho_2$ depends on the importance she attaches to each contrast.
[See Remark 3.3].

Case 2. $ab$ divides $n$, so that orthogonality between $A$ and $B$
is possible. However, in the only available plan (say $\rho_1$) in
which $A$ and $B$ are mutually orthogonal, various other pairs of
factors are mutually non-orthogonal. Suppose another plan $\rho_2$
is also available in which $A$ is only partially orthogonal to $B$,
but several pairs of factors which are mutually non-orthogonal in
$\rho_1$ are orthogonal in $\rho_2$. Which plan should the
experimenter choose ? Again, The choice depends on the importance
attached to different contrasts of different factors. [See  Remark
5.1].

 We hope that in future more and more nearly orthogonal, inter-class orthogonal
and other similar plans will be available and the experimenters will
have a wider range of options.

\vspace{.5em}

 {\bf  Orthogonality through another factor:}


 The concept of ``orthogonality (between two treatment factors)
through a nuisance factor" has been introduced in Morgan and Uddin
(1996) in the context of nested row-column designs. In Bagchi (2010)
``orthogonality through the block factor (OTB)" is studied in
details. This concept can easily be extended to the case when the
third factor is also a treatment factor.

\begin{defi} \label{orth3rd} Consider three factors $A$, $B$ and $C$ of an
MEP. We say that {\bf $A$ is orthogonal to $B$ ``through" $C$} if
the incidence matrices $N_{AB}, N_{BC}$ and $N_{AC}$ satisfy the
following condition.

\begin{equation}  \label{orththr-inc}
 N_{AC} (R_C)^{-1} N_{CB} =   N_{AB}.
\end{equation} \end{defi}

{\bf Example 2.1:} Consider two MEPs with two-  and three-level
factors,  on 5 runs : $A_5(1) = \rho(5; \{3 \times 2^2 \})$ and
$A_5(2) = \rho(5,3;\{2^4\})$.

\begin{equation}  \label{ORTHthruD}
A_5(1) = \left [\begin{array}{cccccc}
A& 0 & 1 &0  & 1 & 0\\
B& 0 & 1 & 1 & 0 & 0\\
C& 0 & 0 & 1 & 1 & 2\end{array} \right ] \;\; A_5(2) = \left
[\begin{array}{cccccc}
 A& 0 & 0 & 1 & 1 & 0\\
 B& 0 & 1 &0  & 1 & 0\\
 C& 0 & 1 & 1 & 0 & 0 \\
 D& 0 & 0 & 0 & 0 & 1
\end{array} \right ].\end{equation}

In the plan $A_5(1)$ $A$ is orthogonal to $B$ through $C$, while in
$A_5(2)$ every factor in $\{A,B,C\}$ is orthogonal to every other
through $D$.

 For examples of more such plans, see the equations next to (\ref {3^2.2^2}).
For analysis see Theorem \ref {POFm}. See also Remark 4.5.

\section{Construction of inter-class orthogonal plans}

\begin{defi} \label{replc-ary} Consider an MEP {\boldmath $\rho(n,m; a,b, \cdots)$}.
Suppose there exists another MEP
 {\boldmath$\rho_1(a,l;t_1, t_2, \cdots t_l)$}, ($l \geq 2$, such that
\begin{equation} \label{sum-lvl}
\sum_{i=1}^{l} (t_{i} -1) \leq a -1. \end{equation}

 Then, we construct a new MEP  {\boldmath $\tilde{\rho}$} with $n$ runs by
replacing the level $u$ of factor $A$  by the $u$-th column (run) of
$\rho_1$, for each $u$, $0 \leq u \leq a-1$. We say that {\bf the
factor A is replaced by a class $G_A$ of $l$ factors related through
$\rho_1$ and $\rho_1$ will be said to the replacing array for $A$}.
In the same way we can replace two or more factors of a given MEP,
through two or more suitable replacing arrays.
\end{defi}

We now try to find conditions on the replacing array so that the
resultant plan  satisfies certain desirable properties.

\begin{lem} \label{replacing}Consider a set of factors ${\cal R} = \{A,B, \cdots\}$ of an MEP $\rho$.
Suppose $\tilde{\rho}$ is an MEP obtained from $\rho$ by replacing
each factor of ${\cal R}$ by a group of factors. More precisely, the
factor A (respectively B ) is  replaced by the class of factors
$G_A$ (respectively $G_B$)  related through $\rho_A$ (respectively
$\rho_B$). Then, the following hold.

(a) If A and B are mutually orthogonal (with equal or unequal
frequency ) in the original plan $\rho$ then every factor in the
class $G_A$ is orthogonal to every factor in the class $G_B$ in the
derived plan $\tilde{\rho}$, generally with unequal frequency.

(b) In $\tilde{\rho}$ two factors of $G_A$ will be partially
(respectively totally) orthogonal  if and only if the corresponding
rows of $\rho_A$ are partially (respectively totally) orthogonal.

(c) If $\rho$ and each of the replacing arrays $\rho_A$ etc. are
saturated, then so is $\tilde{\rho}$.
\end{lem}

{\bf Proof :} We shall prove (a). (b) will follow by similar
argument and (c) by straightforward counting.

 Proof of (a): Fix a factor, say $K$ of $G_A$ and a factor $L$ of $G_B$.
Let $\beta_s$  (respectively $\gamma_t$) denote the set of runs of
$\rho_A$ (respectively $\rho_B$) in which the level s
 of $K$ (respectively t of $L$) appear.

  Let $\tilde{N}_{K,L}, \tilde{r}_K, \tilde{r}_L$ denote the incidence
matrix of $K,L$ and the replication vectors of $K$ and $L$
respectively in $\tilde{\rho}$. Then, the $(s,t)$th entry of
$\tilde{N}_{K,L}$ is given by
\begin{equation} \label{incmatKL} {\tilde{n}_{K,L}} (s,t) =
\sum \limits_{i \in \beta_s} \sum \limits_{j \in \gamma_t} {n_{A,B}}
(i,j). \end{equation}

From this, we obtain that for a level $s$ of $K$,
\begin{equation}
\tilde{r}_K (s) \;  = \; \sum \limits_{i \in \beta_s} \sum
\limits_{t}
 \sum \limits_{j \in \gamma_t} n_{A,B} (i,j) \;
     =\;  \sum \limits_{i \in \beta_s} r_A (i). \end{equation}
Similarly, $\tilde{r}_L (t) = \sum \limits_{j \in \gamma_t} r_B
(j)$. But $N_{A,B},r_A,r_B$ satisfy (\ref {PFC}) by hypothesis.
Combining that with the  relations above we see that
$\tilde{M}_{K,L}, \tilde{r}_K, \tilde{r}_L$ also satisfy (\ref
{PFC}). $\Box$

\vspace{.5em} We present the well-known definition of an orthogonal
array.

\begin{defi} \label {OAnmst}Let $m,n,t \geq 2$ be integers and $s = (s_1, \cdots
s_m)$ be a vector of integers $\geq 2$. Then an orthogonal array of
strength $t$  is an $m \times n$ array, with the entries of the
$i$th row coming from a set of $s_i$ symbols satisfying the
following. All $t$-tuples of symbols appear equally often as rows in
every $n\times t$ subarray.
 Such an array is denoted by $OA(m,n, s_1\times \cdots \times s_m,
t)$. When $ s_1 = s_2 = \cdots s_m =s$, say, this array is
represented by $OA(n,m,s,t)$.
\end{defi}

\begin{cor} \label{rplc-OA} Suppose there exists an orthogonal
array $OA(n,m,s,2)$. Suppose further for an integer $k (<m)$, there
exist  arrays $\rho_i =\rho(s,l_i;t_{i,1}, \cdots t_{i,l_i})$,
satisfying $s-1 \geq \sum_{j=1}^{l_i} (t_{i,j} -1), i=1,2, \cdots
k$.

 Then, an inter-class orthogonal array $\rho(n,l;  \prod_{i=1}^k
\{t_{i,1} \times \cdots t_{i,l_i}\} . s^{m-k}) $ exists. Here
 $l = \sum_{i=1}^{k} l_{i}$.
 \end{cor}

 {\bf Proof :} Let $\rho_0$ be the orthogonal MEP represented by the
 given orthogonal array. We replace the ith factor by a group of
 factors related through $\rho_i, \: i=1, \cdots k$ to form a new
 MRP  $\rho$. Clearly  $\rho$ has $l = \sum_{i=1}^{k} l_{i} + m - k$
 factors. That $\rho$ is interclass orthogonal with the given
 parameters follows from Lemma \ref {replacing}.$\Box$

\vspace{.5em}

{\bf Examples of replacing arrays with desirable properties: }

\vspace{.5em}

 We have seen that to obtain an useful inter-class orthogonal MEP,
one needs replacing arrays with desirable properties. We now present
a few such arrays. In each plan, the factors are named as $A, B,
\cdots$, in that order. The set of $s$ levels of a factor will be
denoted by the set of integers modulo $s$.

{\bf Plan with $s$ runs, two factors with $p$ and $q$ levels,
$p+q=s$:}

\begin{equation}\label{genRP} A_s(1) = \left [ \begin{array}{cccccccc}
                  0 & 0 & \cdots & 0 & 1 & 2 & \cdots & p-1\\
                  0 & 1 & \cdots & q-1 & 0 & 0 & \cdots & 0 \\
                           \end{array} \right ]. \end{equation}

{\bf Remark  3.1:} The BLUE of the contrast $\alpha_i - \alpha_j$ of
factor $A$ is orthogonal to the BLUEs of the contrasts of B, for $i
\neq j, \: i,j \geq 1$. Similarly, the BLUE of the contrast $\beta_i
- \beta_j$ of $B$ is orthogonal to the BLUEs of the contrasts of A,
for $i \neq j, \: i,j \geq 1$.

\vspace{.5em}

{\bf Plans with 4 runs :} A plan, say $A_4(1)$ may be obtained by
putting $s = 4, p = 2, q = 3$ in (\ref {genRP}). We now present
another plan.

$A_4(2) = \rho(4; \{3 \times 2\})  = \left [\begin{array}{cccc}
 0 & 1 & 0 & 2 \\
  0 & 1 & 1 & 0 \\\end{array} \right ]$.

\vspace{.5em}

{\bf Remark 3.2 :} (a) In $A_4(2)$, (\ref {PFCi}) is satisfied by
level 0 of factor $A$, so that  $2\hat{\alpha} - \hat{\alpha_1} -
\hat{\alpha_2}$ is orthogonal to the BLUEs of the main effect
contrast of $B$.

\vspace{.5em}

{\bf Plans with 5 runs :} Two plans, namely $A_5(1) = \rho(5; \{2^2
\times 3\})$ and $A_5(2) = \rho(5; \{2^4\})$ are presented in
Example 2.1. Two other plans $A_5(3)$ and $A_5(4)$ are obtained from
(\ref {genRP}) by putting $s = 5, p = 4, q = 2$ and $s=5, p= q =3$
respectively. \vspace{.5em}

{\bf  Plans with 7 runs :}  A plan, say $A_7(1)$ may be obtained by
putting $s = 7, p = 6, q = 2$ in (\ref {genRP}). Another plan is
displayed below.

$A_7(2) = \rho(7; \{3^3\})  = \left [\begin{array}{ccccccc}
 0 & 0 & 0 & 1 & 1 & 2 & 2 \\
 0 & 1 & 2 & 0 & 1 & 0 & 2 \\
 0 & 1 & 2 & 1 & 2 & 2 & 0 \end{array} \right ]$.

\vspace{.5em}

Another plan, say $A_7(3) = \rho(7; \{4^2\})$ may be obtained by
taking the first two rows and the columns numbered 2,3,7,8,10,12,and
13 from the array $R^7$ in (3.17).

 In the next section we use suitable arrays from the list above
 to replace one or more rows of existing orthogonal arrays and obtain inter-class
orthogonal MEPs. Before that we compare  the two replacing arrays
$A_4(1)$ and $A_4(2)$ regarding the precision of the BLUEs of the
main effect contrasts. We first compute the C-matrices (the
coefficient matrices). $C_{AA;\bar{A}}$ denotes the coefficient
matrix of the system of reduced normal equations for factor $A$.
[See Notation \ref {C,QandSS} and (c) of Corollary \ref {rdcedNE}].

It is rather surprising that $C_{BB;\bar{B}}$ is the same for both
the plans. $ C_{BB;\bar{B}} = \left [  \begin{array} {ll} 1/2 & -1/2
\\ -1/2 & 1/2 \end{array} \right ].$

$ C_{AA;\bar{A}}$ is, however different in the two plans. They are
given below.

$$C_{AA;\bar{A}} \mbox{ for $A_4(1)$ is }
 \left [ \begin{array}{ccc}2/3 &  -1/3 & -1/3\\
                           -1/3 &  2/3 &  -1/3 \\
                           -1/3 & -1/3 & 2/3 \\
\end{array} \right ] \mbox{ and  $C_{AA;\bar{A}}$ for $A_4(2)$ is }
 \left [ \begin{array}{ccc}
1      &  -1/2 & -1/2 \\
 -1/2  &  1/2  & 0      \\
 -1/2 &  0     & 1/2`. \\
\end{array} \right ].$$

{\bf Remark 3.3:} We note that both contrasts of the three-level
factor $A$ are estimated with the same precision in $A_4(1)$. In
$A_4(2)$, however, the contrast $2\hat{\alpha}_0 - \hat{\alpha}_1 -
\hat{\alpha}_2$ orthogonal to the BLUEs of the contrasts of $B$ and
hence is estimated with maximum possible precision (given the
replication vector), while the contrast $ \hat{\alpha}_1 -
\hat{\alpha}_2$  is estimated with much less precision. Thus, while
replacing a four-level factor the experimenter may choose between
$A_4(1)$ and $A_4(2)$ depending on whether equal importance is
attached to both the contrasts or not.

\vspace{.5em}

\subsection{ Some series of inter-class orthogonal main effect plans}
Our starting point is an {\it $OA(n,m,s,t)$} (see Definition \ref
 {OAnmst}).

\begin{theo}\label{OAgen}
(a) Whenever an $OA(n,m,\prod_{i=1}^{m} s_{i},2)$ exists, an
inter-class orthogonal MEP $\rho(n,2m; \prod_{i=1}^{m} \{(s_i
-t_i).(t_i +1) \}$ exists. Here $t_i$ is an integer,  $1 \leq t_i
\leq s_i-2$.

 (b) These inter-class orthogonal MEPs may be constructed so as to
 satisfy partial orthogonality property among the members
 of the same class similar to the description in Remark 3.1.
\end{theo}

{\bf Proof :} (a) For every $i, 1\leq i \leq m$, one can choose a $
2 \times s_i$ array, say $\rho_i$, with $p = s_i - t_i$ symbols in
the first and $q = t_i + 1$ symbols in the second row.
 Now  $\rho_i$ may be used as a replacing array for the
ith factor of the given OA.

(b) In particular, if  $\rho_i$ has the same structure as $A_s(1)$,
with $s=s_i, p = s - t_i, q = t_i + 1$, then the members of the ith
class will satisfy the stated partial orthogonality property.

\begin{theo} \label{OA(n,m;s^m)} Suppose $s=3,4,5$ or 7. Whenever an
$OA(n,m,s,2)$ exists, the following  series of inter-class
orthogonal MEPs $\rho_1$ exist. Here $p,q,r,s,t$ are nonnegative
integers.

\begin{equation} \label{1st-array}
 \rho_1 = \left\{ \begin{array}{lll} \rho(n; 3^p \times \{2^2\}^q),
 & p+q=m,   & \mbox{ if s=3} \\
 \rho(n; 4^p.\{3.2\}^q.2^{3t}) & p+q+t=m, & \mbox{ if s=4} \\
 \rho(n; 5^p.\{4.2\}^q. \{3^2 \}^r.\{3.2^2\}^t.\{2^4 \}^u) &
  p+q+r+t+u=m,  & \mbox{ if s=5} \\
 \rho(n; 7^p.\{6.2\}^q. \{4^2 \}^r.\{3^3\}^t) &
  p+q+r+t=m,  & \mbox{ if s=7}
\end{array} \right. \end{equation}\end{theo}

{\bf Proof: } Let $O$ be an $OA(n,m,s,2)$. We keep $p$ (out of m) of
the factors of $O$ as they are, replace every other factor by a
class of factors related through an appropriate replacing array.
This replacing array can be (i) an OA if it is available (which is
the case when $s=4$), (ii) one of the replacing arrays shown above
or (iii) a replacing array of similar type, for instance $A_3(1)$,
obtained by putting $s=3$ in $A_s(1)$. Corollary \ref {rplc-OA}
implies that the MEP thus constructed satisfies the required
property. $\Box$

{\bf Discussion :} 1. While applying Theorem \ref {OA(n,m;s^m)} with
$s=4$, the experimenter has a choice between the replacing arrays
$A_4(1)$ and $A_4(2)$. Remark 3.3  may be useful in making the
choice.

2. Comparing an inter-class(2) orthogonal MEP, say  $\rho_1$
constructed in Theorem \ref {OA(n,m;s^m)} with $s=4$ with an
existing plan with the same number of runs, we find the following.
In the plan $ICA(n,3^l,2^{n-2-2l})$ of Huwang, Wu and Yen (2002),
 a three-level factor is orthogonal to every two-level factor and
 non-orthogonal to every other three-level factor [see p-349,
 line 7 of HWY]. In  $\rho_1$ every three-level factor is orthogonal
 to every other three-level factor
 and all but one two-level factors, (with which it is partially
 orthogonal in case $A_4(2)$ is used).

\subsection{More series of inter-class (2) orthogonal main effect plans}

We shall now present a two-stage construction. In the first stage we
start with an existing MEP, fix a subset (say ${\cal R}$) of factors
and obtain a number of MEPs by replacing  each factor in  ${\cal R}$
by a class of factors. Here we may use  different replacing arrays
for the same factor while constructing the first stage MEPs.

In the next stage we juxtapose these the first stage MEPs in a
suitable manner to form an array. In order that the resultant array
is a meaningful MEP, the replacing arrays need to satisfy certain
condition as we shall see now.

\begin{defi} Consider an MEP $\rho$. Let ${\cal F}$ denote the class
of all factors of $\rho$. Suppose $\rho_P(1)$ and $\rho_P(2)$ denote
two replacing arrays for a factor $P$. If these replacing arrays
have the same number of factors  and the same number of levels for
the corresponding factors, then they are said to be compatible.
Further both of them are said to represent the same class of
factors, say $G_P$.

Let $\rho_1$ and $\rho_2$ be two MEPs obtained from $\rho$ by
replacing the factors in a certain subset ${\cal R}$ of ${\cal F}$.
$\rho_1$ and $\rho_2$ are said to be compatible w.r.t. the factor
$P$, if the corresponding replacing arrays $\rho_P(1)$ and
$\rho_P(2)$ for $P$ are compatible, in which case we say that both
$\rho_1$ and $\rho_2$ are obtained by replacing $P$ with the same
class $G_P$ of factors.

If $\rho_1$ and $\rho_2$ are compatible w.r.t. each factor in ${\cal
R}$ then they are said to be compatible w.r.t ${\cal R}$. \end{defi}

The following results are immediate from the definition.

\begin{lem} \label{rplc-jn1}  Consider an MEP $\rho$ and a subset
 ${\cal R}$ of ${\cal F}$. Suppose for every factor $P$ in
 ${\cal R}$ of $\rho$ there is a class of
replacing arrays $\rho_P(i,j), j=1,2, \cdots J, i=1,2, \cdots I$,
such that the replacing arrays $\rho_P(i,j), j=1,2, \cdots J$ are
mutually compatible for every $i=1,2, \cdots I$ and every $P$ in
${\cal R}$. For $j=1,2, \cdots J$ let $G_P(i)$ denote the class of
factors of each $\rho_P(i,j), j=1,2, \cdots J$.

Now, for every $i=1,2, \cdots I, j=1,2, \cdots J$, we obtain an
array $\rho_{ij}$ by replacing a factor $P$ of $\rho$ by a class of
factors $G_P (i)$ related through $\rho_P(ij)$. Let

\begin{equation} \label{newArray} \rho^* = ((\rho_{ij}))_{1 \leq i
\leq I,\: 1 \leq j \leq J}.\end{equation}

Then $\rho^*$ represents an MEP satisfying the following.

(a) $\rho^*$ can be viewed as an MEP directly obtained from $\rho$
by replacing every $P$ in ${\cal R}$ with a class $G_P$ of
 factors related through the following array.

\begin{equation} \label{fullRplcArr} \rho^*_P = ((\rho_P(ij)))_{1 \leq i
\leq I,\: 1 \leq j \leq J}.\end{equation}

(b) If P and Q are mutually orthogonal in the original plan $\rho$ ,
then every factor in the class $G_P (i)$ is orthogonal to every
factor in the class $G_Q (i')$ in the derived plan $\rho^*$, $
i,i'=1,2, \cdots I$.

(c) Fix a factor $P$ of $\rho$. Fix $i\neq i',\; i,i' =1,2, \cdots
I$. let

$$\rho_P (i,i') =\left [ \begin{array}{ccc}
      \rho_P (i1) & \cdots & \rho_P (i,J) \\
      \rho_P (i'1)& \cdots & \rho_P (i'J) \\
         \end{array} \right ].$$

Let us look at the set of factors $G_P (i) \cup G_P (i')$ of  the
derived plan $\rho^*$. Two factors  in this class are mutually
orthogonal if and only if the corresponding factors in the plan
represented by $\rho_P (i,i')$ are so. \end{lem}

{\bf Remark 3.4:} In our definition of replacing arrays we have used
the condition (\ref {sum-lvl}), so that they are not supersaturated
and hence the resultant MEPs are also not supersaturated. However,
we now relax this condition a bit. That is, we make use of one or
more supersaturated replacing arrays in the intermediate stage, but
the final MEP will not be supersaturated.

\begin{lem} \label{rplc-jn2} Consider a set up just like that in
the statement of Lemma \ref {rplc-jn1}, except the following. There
is a factor $P$ and an i, say $i_0$, such that $\rho_P(i_0j)$ is
supersaturated, i.e. it does not satisfy (\ref {sum-lvl}) for every
$j = 1, \cdots J$.

  Let $\rho^*$ and $\rho_P^*$ be  as in Lemma \ref {rplc-jn1}.
Let $ G_P$ denote the class of factors $ G_P = \bigcup_{i=1}^I
G_P(i)$.

   Then, statements (a) and (b) of Lemma \ref {rplc-jn1} hold.
Further, the following modified form of Statement (c) of the same
Lemma hold.

(c)'  If $\rho^*_P$ [see (\ref {fullRplcArr})] is not
supersaturated, then

(i) in $\rho^*$ any pair of  factors in $ G_P$ are mutually
orthogonal if and only if the corresponding factors in the plan
represented by $\rho^*_P$ are so and

(ii)  the main effect contrast of each member of $ G_P$ can be
estimated.
\end{lem}

We now apply the technique of two stage construction to construct
more inter-class orthogonal MEPs with two or three levels. Some of
them turn out to be (fully) orthogonal.

\begin{theo} \label{two-stage} (a) The existence of an $OA(n,m,s,2),\;
s=4,5$ or 7 implies the existence of the following inter-class
orthogonal MEP.

\begin{equation} \label{2nd-array}
 \rho_2 = \left\{ \begin{array}{llll}
\rho(2n,4m; \{3^2\}^m.2^{2m})& & &\mbox{ if s=4} \\
\rho(2n,4m; \{3^4\}^m) & & &\mbox{ if s=5} \\
\rho(2n,4m; \{4^4\}^m) & & &\mbox{ if s=7}
\end{array} \right.\end{equation}

Further, this MEP satisfies the following properties.

(b)  In the case $s=4$, the pairs of three-level factors are
 partial orthogonal to each other - in fact the relation between
the pairs of three-level factors is just like the factors $A$ and
$B$ of $A_8(1)$ [See Section 1].

so that every contrast is orthogonal to all except possibly another
contrast. In particular, the relation between the

(c) In the cases $s=5$ and $s=7$, among the four members in the same
class, every pair among the last 3 are mutually orthogonal through
the first one. \end{theo}

{\bf Proof :} Fix $s \in \{4,5,7\}$ Let $R^s(4 \times 2s)$ denote a
suitably chosen array, which is partitioned as

$$ R^s = ((R_{ij}))_{1 \leq i,j \leq 2}, \mbox{ each  $R_{ij}$ is } 2 \times s.$$

 Let $O$ denote the given OA. We first construct four arrays
$\rho_{11}, \rho_{12}, \rho_{21}$ and $ \rho_{22}$ following the
method of Theorem \ref {OA(n,m;s^m)}. In this process we use
$R_{ij}$ as the replacing array for each factor to construct
$\rho_{ij}, i,j =1,2$. Now we form

$$ \rho^*  = ((\rho_{ij}))_{1 \leq i,j \leq 2},$$

which is the required MEP. By Lemma \ref {rplc-jn1}, it follows that
$ \rho^*$ may be viewed as the plan obtained by replacing every
factor P by the class $G_P$ of four factors related through the
replacing array $R^s$. The rest of the proof follows from the
structures of $R^s, s=4,5,7$, shown below.

\begin{eqnarray} \label{3^2.2^2}
 R^4 &= &\left [\begin{array}{cccccccc}
    0 & 0 & 1 & 0 & 0 & 1 & 1 & 1 \\
    0 & 1 & 0 & 2 & 0 & 1 & 0 & 2 \\
    1 & 0 & 0 & 1 & 0 & 1 & 1 & 0 \\
    2 & 1 & 0 & 0 & 0 & 0 & 1 & 2   \end{array}\right ]\\
R^5 &= &\left [\begin{array}{cccccccccc}
0  & 0 &  0 &  0 &  1 &  1 &  1 &  1 &  2 &  2\\
0  & 0 &  1 &  2 &  0 &  0 &  1 &  2 &  1 &  2 \\
0  & 1 &  0 &  2 &  2 &  0 &  1 &  0 &  2 &  1\\
0  & 1 &  2 &  0 &  2 &  0 &  0 &  1 &  1 &  2\end{array}\right
]\\
R^7 &= &\left [\begin{array}{cccccccccccccc}
0  & 0 &  0 &  0 &  1 &  1 &  1 &  1 &  2 &  2 & 2 & 2 & 3 & 3\\
0  & 0 &  1 &  2 &  0 &  0 &  1 &  2 &  1 & 1 &  3 & 3 & 1 &  2 \\
0  & 1 &  0 &  2 &  2 &  0 &  1 &  0 &  1 & 3 & 1 & 3  & 2 &  1\\
0  & 1 &  2 &  0 &  2 &  0 &  0 &  1 &  1 & 3 & 3 & 1  & 1 &2
\end{array}\right ] \Box \end{eqnarray}

   Our next result is based on the elegant plan of Stark (1964), which
   is quoted below. [See Dey (1985), for instance, for an explicit
   presentation of the plan and more details].

\begin{theo} (Stark (1964))  An OMEP for a $3^7$ experiment on 16 runs
exists. \end{theo}

\begin{theo} \label{two-stage-Stark} (a)The existence of an $OA(n,m,8,2)$
 implies the existence of an orthogonal MEP for $7m$ three-level
 factors on $2n$ runs.

(b) The existence of an $OA(n,m,4,2)$ implies the existence of an
orthogonal MEP for $6m$ three-level factors on $4n$ runs.
\end{theo}

 {\bf Proof:} Let $R$ be a $7 \times 16$ array with symbols $0,1,2$
 representing the OMEP of Stark.

 (a) We partition $R$ as $ R =\left [ \begin{array}{cc}
    R_1 & R_2 \end{array} \right ], \; \mbox{ each } R_i \mbox{ is of order }
7 \times 8 .$

  We first construct arrays $\rho_{j}$ from the given
 OA by using replacing array $R_{j}$ for every factor following the
 method of Theorem \ref {OA(n,m;s^m)} $:j=1,2$. Then we form the required plan $\rho^*$
as
$$\rho^* = \left [\begin{array}{cc}\rho_1 & \rho_2\end{array} \right
].$$

  That $\rho^*$ satisfies the required property follows from Lemma
  \ref {rplc-jn2}.

(b) Let $\tilde{R}$ denote the $6 \times 16$ array obtained by
deleting a row (say the 0th one) from $R$. Now we partition
$\tilde{R}$ as follows.

$$\tilde{R} =  ((\tilde{R}_{ij}))_{1 \leq i,j \leq 4},$$

such that $\tilde{R}_{ij}$ is of order 2 $\times$ 4 for $i=1,2$ and
1 $\times$ 4 for $i=3,4$.

 Let $A$ denote the given OA. We  construct array $\rho_{ij}$
 from $A$ by using replacing array $\tilde{R}_{ij}$ for every factor
 following the method of Theorem \ref {OA(n,m;s^m)}, $i,j=1,2,3,4$. Then we form
the required plan $\rho^*$ as

 $$ \rho^* = ((\rho_{ij}))_{1 \leq i,j \leq 4}.$$

 Note that this procedure may be viewed as follows. Fix a factor,
 say $P$ of $A$. The intermediate arrays $\rho_{ij}, 1\leq j \leq 4, i=1,2$
 are formed by replacing every factor $P$ of $A$ by two three-level factors
 each, so that $\rho_{ij}, 1\leq j \leq 4, i=1,2$ are supersaturated. However,
 in each of the  intermediate arrays
 $\rho_{ij}, 1\leq j \leq 4, i=3,4$, $P$ is replaced by one
three-level factor. This fact,
 together with the choice of replacing arrays imply that the class
 $G_P$ (the class of factors in $\rho^*$ replacing P) is nothing
 but a class of six three-level factors, related through $\tilde{R}$.
 The rest follows from Lemma \ref {rplc-jn2} and the fact that
 there are $m$ factors in $A$. $\Box$

\section{Analysis of a general main effect plan.}

 The crucial component of data analysis of a general factorial experiment
  is, of course, the computation of the error sum of squares. We
  proceed towards a user-friendly formula for computing $SS_E$. The results
  are not new, but are not available in the form presented here.
 We denote the factors by $F_1, F_2, \cdots $,
     instead of $A,B \cdots$ for the sake of notational simplicity.

 We assume an additive, fixed effects, main effects model with
homoscedastic and uncorrelated errors having constant variance
$\sigma^2$. ${\mathbf 1}_n$ will denote the $n \times 1$ vector of
all-ones, while ${\mathbf J}_{m \times n }$ will denote the $m
\times n$ matrix of all-ones.

 Let $\rho$ denote a main effect plan on $n$ runs with  factors $F_1, F_2, \cdots F_m$,
 $F_i$ having $a_i$ levels, $i=1, \cdots m$. Let the  unknown effect of
 the jth level of the factor $F_i$ be denoted by $\alpha^i_j$ and let the
 $ a_i \times 1$ vector $\alpha^i$ denote the vector of unknown
 effects of $F_i,\: 1\leq i \leq m$.
 Let  $Y_u$ denote the yield from the
 uth run, $u=1,2,\cdots n$. Then, assuming that in the $u$th run
 the factor $F_i$ is set at level $l_i = l_i(u), i=1, \cdots m$ and denoting
 the general effect by $\mu$, $Y_u$ is given by
$$ Y_u = \mu  + \sum_{i=1}^{m} \alpha^i_{l_i} + \epsilon_u,
u=1,2, \cdots n. $$

Viewing the general effect as the (m+1)-th factor ( $F_{m+1}$)
 and therefore writing $\alpha^{m+1} = \mu$ we express the model in
matrix form as
\begin{equation}\label{model} {\mathbf Y} =  {\mathbf X} \beta,
\mbox{ where } {\mathbf X}  = \left[ \begin{array}{ccc}  {\mathbf
X}_1 & \cdots &  {\mathbf X}_{m+1}\end{array} \right ]  \mbox{ and
}\beta =\left[ \begin{array}{ccc} \alpha^1 & \cdots &
\alpha^{m+1}\end{array} \right ]^T.\end{equation}

Here, $X_i$, the {\bf design matrix} for $F_i$ is a $0-1$ matrix -
the $(u,t)$th entry of ${\mathbf X}_i $ is 1 if in the uth run the
factor $F_i$ is set at level t and 0 otherwise, $ i=1,2, \cdots m$
and $X_{m+1} = {\mathbf 1}_n$.

Let  ${\mathbf T}_i$ denote the vector of raw totals of $F_i,\: i=1,
\cdots m+1$. Thus, ${\mathbf T}_{m+1}$ is the grand total and will
sometimes be denoted by $G$.

\begin{nota}  For any $m \times n$ matrix $A$, ${\cal C} (A)$
will denote the column space of $A$. Further,  $P_A$ will denote the
projection operator on the column space of $A$. In other words, $P_A
= A (A'A)^- A'$, where $B^-$ denotes a g-inverse of $B$.
\end{nota}

\begin{nota}\label{setofFact}
  Let $I =  \{1,2, \cdots m +1\}$ and   $S = \{i,j, \cdots\}$ be
  a subset of $I$. For the sake of compactness, we introduce the following
 notation.

 (a) $\bar{i} = I \setminus \{i \}$.

 (b) ${\mathbf X}_S = \left[ \begin{array}{ccc}{\mathbf X}_i &
 {\mathbf X}_j & \cdots \end{array} \right ]$.

(b) $ \alpha^S = \left[ \begin{array}{ccc} \alpha^i & \cdots &
\alpha^j  \end{array} \right ]^T$.

(c)  ${\mathbf P}_i$ will denote the projection operator onto the
column space of ${\mathbf X}_i, i \in I$. Further, ${\mathbf P}_S$
will denote the projection operator onto the column space of
${\mathbf X}_S$.
\end{nota}

 {\bf The system of reduced normal equations for a class of factors.}

\begin{nota} \label{C,QandSS}(a) Let $S,T,U$ be three subsets of $I$
such that

(i) $S \cap U = T \cap U = \phi$ and

(ii) either $S = T$ or $ S \cap T = \phi$.

 Let us define the matrix $C_{S,T;U}$ and
the vector $Q_{S;U}$ as follows.
\begin{eqnarray} \label{CQSS}
 {\mathbf C}_{S,T;U} & = & (({\mathbf C}_{ij;U}))_{i \in S,\: j \in T},\;
 {\mathbf C}_{ij;U}  =  {\mathbf X}'_i (I - P_U) {\mathbf X}_j, \\
{\mathbf Q}_{S;U} & = & (({\mathbf Q}_{i;U} ))_{i \in S}, \;
{\mathbf Q}_{i;U} = {\mathbf X}'_i (I - P_U) {\mathbf Y}.
\end{eqnarray}

 (b) In particular, if S and T are a singleton sets, say $S = \{i  \}$
and $T = \{j \}$ , then we may write ${\mathbf C}_{ij;U}$ and
${\mathbf Q}_{i;U}$ instead of $ {\mathbf C}_{S,T;U}$ and ${\mathbf
Q}_{S;U}$ respectively. Sometimes we may write $C_{i;U}$ instead of
$C_{ii;U}$.

(c)  Suppose $U = \{k,m+1\}$. Then we may and do write ${\mathbf
C}_{S,T;k}$ and ${\mathbf Q}_{S;k}$ instead of ${\mathbf C}_{S,T;U}$
and ${\mathbf Q}_{S;U}$ respectively. [This is because ${\mathbf
P}_U = {\mathbf P}_k$ ].

\end{nota}

  The following  well-known result is presented using the notation above.

\begin{lem} \label{CandQ}  Suppose $I$ is partitioned into two
subsets  $S$  and $U$. Then, the reduced normal equations for
$\widehat{\alpha^S}$, after eliminating $\widehat{\alpha^U}$ is
given by
$$ C_{S,S;U} \widehat{\alpha^S} = Q_{S;U},$$

where $ C_{S,S;U}$  and $ Q_{S;U}$ are as given in ( \ref{CQSS}) and
the next equation.
\end{lem}

  {\bf Remark 4.1:} In order that every  main effect contrast of $F_i$
is estimable, rank of $ C_{i;\bar{i}}$ must be $a_i - 1$, as we
know. Thus, before using a general m-factor MEP $\rho$ with $m \geq
3$, one has to check whether  $ Rank (C_{i;\bar{i}}) = a_i - 1$, for
every $i =1,2, \cdots m$.

 In view of Remark 4.1 above, we define a class of MEPs,
 borrowing a term from the theory of block designs.

\begin{defi}\label{connected} An m-factor MEP is said to be
`connected' if $Rank( C_{i;\bar{i}}) = a_i - 1$, for every $i =1,2,
\cdots m$.
\end{defi}

   Henceforth, the MEP $\rho$ under consideration will be assumed to
   be connected. We now present a few special cases of Lemma \ref {CandQ}

\begin{cor} \label{rdcedNE} (a) Consider a factor, say $F_i$. Let

 $\bar{i} = I
\setminus \{i \}$. Then the BLUE of the main effect contrast
 $l'\alpha^i$ (in case it is estimable) of $F_i$
is $ l'\widehat{\alpha^i}$, where $\widehat{\alpha^i}$ is a solution
of
  $$ C_{i;\bar{i}} \widehat{\alpha^i} = Q_{i;\bar{i}}. $$

Here the expressions for $ C_{i;\bar{i}}$ and $ Q_{i;\bar{i}}$ are
obtained from (b) of Notation \ref {C,QandSS}.

 (b) In particular, suppose $m = 1$. Then, the reduced normal equation
 for $\alpha^1$ (obtained by eliminating only $F_2 =\mu$) is
\begin{equation}\label{C01}
({\mathbf R}_1 -  {\mathbf r}_1  ({\mathbf r}_1)'/n)
\widehat{\alpha^1} = {\mathbf T}_1 - {\mathbf r}_1 G/n
.\end{equation}

(c) Suppose $m = 2$. Then the reduced normal equation for $\alpha^1$
(obtained by eliminating $F_2$ and  $F_3 = \mu$) is
$C_{1;2}\widehat{\alpha^1} = Q_{1;2}$, where
\begin{equation} \label{CQ12}
  C_{1;2}= {\mathbf R}_1 - N_{12} (R_2)^{-1} N_{21} \mbox{ and }
  Q_{1;2} = {\mathbf T}_1 - N_{12}( R_2)^{-1} {\mathbf T}_2.
\end{equation}
\end{cor}

\begin{nota}\label{adjSS} We now define sum of squares for one or more factors,
adjusted for one or more other factors. Fix a set of factors $T$ of
$I$. For i not in $T$, we define $SS_{i;T}$, the sum of squares for
$F_i$, adjusted for the factors $F_t, t \in T$. More generally, for
S disjoint from T, we define  $SS_{S;T}$, the sum of squares for the
set of  factors $F_i, i \in S$, viewed as a single factor, adjusted
for the factors $F_t, t \in T$.

\begin{eqnarray*} SS_{i;T} & = & Q'_{i;T} (C_{i;T})^- Q_{i;T} \\
\mbox{ and }  SS_{S;T} & = &Q'_{S;T} (C_{S,S;T})^- Q_{S;T}
\end{eqnarray*}
 \end{nota}

 {\bf Remark 4.2:} Consider two sets of disjoint factors $S$ and $T$.
We may view all the factors in $S$  combined together as a single
factor, say $F_S$,  having design matrix $X_S$. Similarly $F_T$ is
the set of all factors in $T$. Then $SS_{S;T}$ may be viewed as the
sum of squares for $F_S$ adjusted for $F_T$.

 In order to study the relationship between the sums of squares, we
 need the following results on partition matrices.

 \begin{lem} \label{PA-PD} Consider a matrix $W$ partitioned as
$\left[ \begin{array}{cc} U & V \end{array} \right ]$. Let $ Z = (I
- P_V)U$. Then, $P_W - P_V = P_Z$. \end{lem}

\vspace{.5em}

\begin{cor} \label{Ti} Let $ T \subset I$ and $i \in I \setminus T$.
 Let $D = (I - P_T)X_i$. Then,

$$P_D = P_{T*} - P_T, \mbox{ where } T* = T \cup \{i \}.$$
\end{cor}

 We need some more notation.

\begin{nota} \label{sumsofsq}(a) The total sum of squares and the error sum of squares
will be denoted by $SS_{tot}$ and $SS_E$ respectively.

(b) Fix a factor $F_i, 1 \leq i \leq m$. Let $T = \{i+1, \cdots m+1
\}$ and $\bar{i} = I \setminus \{i \}$.

(i) Let $SS_{i;all>} = SS_{i;T}$ and

(ii) $SS_{i; all} = SS_{i;\bar{i} }$.

Thus, $SS_{i;all>}$ is the sum of squares for $F_i$, adjusted for
the factors $ F_{i+1}, \: \cdots F_{m+1}$, while $SS_{i; all}$
denotes the sum of squares for $F_i$, adjusted for all other
factors.
\end{nota}

 {\bf Remark 4.3 :} Note that  $SS_{m;all>}$ is the so-called
 unadjusted  sum of squares for $F_m$.

 We are now in a position to present the computational formulae of
 the error sum of squares.

\begin{theo}\label{FerrorSS} Consider a main effect plan with m
mutually non-orthogonal factors $F_1, F_2, \cdots F_m$.
 The
 error sum of squares ($SS_E$) may be computed from the  total sum of
 squares ($SS_{tot}$) as follows.
\begin{eqnarray} \label{SSE} SS_E    &=& SS_{tot} - SS_{sub},\: \mbox{ where}  \\
                  SS_{sub}& =& \sum_{i=1}^{m} SS_{i;all>}.
\end{eqnarray}

\end{theo}



 \vskip5pt
\begin{theo} \label{ANOVA}  The data obtained from a connected
main effect plan with m mutually  non-orthogonal factors may be
analyzed using the following table.

\begin{center}
{ Table 2.1 : ANOVA for an m-factor non-orthogonal main effect plan}
\end{center}

\begin{tabular}{ccccc}
\hline\\ Source & d.f & S.S adjusted   & S.S adjusted & F-statistics     \\
                &     & for all others & for the next ones &\\
\hline
        &        &           &  &\\
$F_1$ &$a_1-1$ &$SS_{1;all}$ & $SS_{1;all>} = SS_{1;2, \cdots m}$ &
$\frac{SS_{1;all}/(a_1-1)}{SS_E/e}$\\
&  &  & &  \\
$F_2$ &$a_2 - 1$&$SS_{2;all}$         & $SS_{2;all>} =
SS_{2;3,\cdots m}$&
$\frac{SS_{2;all}/(a_2-1)}{SS_E/e}$\\
\vdots & \vdots & \vdots &  \vdots &\\
&  &  &  &\\
$F_m$ & $a_m -1$ & $SS_{m;all}$ & $SS_{m;all>}= SS_{m;m+1}$&
 $\frac{SS_{m;all}/(a_m-1)}{SS_E/e}$ \\
 \hline \\
To be subtracted & -&- & $SS_{sub} =$ Sum of all above\\
&  &  &  &\\
Error &$e$& $SS_E$ =  &$SS_{tot} - SS_{sub}$ \\
&  &  & & \\
Total & n-1 & $SS_{tot}=$ & $\sum_{u=1}^{n} Y^2_u - G^2/n$  & -
\\\hline
\end{tabular}

\vspace{.5em}

Here the error degrees of freedom is $e = n-1-\sum_{i=1}^{m} (a_i -
1)$, as usual. \end{theo}

\vspace{.5em}

{\bf Extension to a general factorial experiment:} Consider a plan
for a factorial experiment with k factors. Let $E$ denote a
factorial effect - a main effect or a t-factor interaction, $ 2 \leq
t \leq k$. We list the factorial effects under study as say $E_1,
\cdots E_m$, where $m$ is the number of factorial effects of
interest. Then we treat these $E_i$'s in the same way as the main
effects $F_i$'s are treated above. That is we denote the design
matrix and the unknown effects of $E_i$ as $X_i$ and $ \alpha^i$ as
here, orders of these would be different when the effects are
interactions. Thus, following the same argument, we arrive at the
following result.

\begin{theo} The error sum of squares of a general factorial
experiment can be obtained in the same manner as described in
Theorem \ref {FerrorSS}.\end{theo}

 {\bf  Situations when analysis is considerably simpler.}

  \vskip5pt

We have seen in Theorem \ref {FerrorSS} that analysis of a general
MEP is rather involved - needs computation of $2m-1$ sums of squares
for an m-factor plan. We now look for situations when so much
computation is not needed.

 We know that when there is only one treatment factor  $F_1$,
the sum of squares for $F_1$ is nothing but the so-called unadjusted
sum of squares $ T'_1 (R_1)^{-1} T_1 -G^2/n$. Moreover, in the
situations when there are two factors, say $F_1$ and $F_2$, the sum
of squares for $F_1$ (obtained by adjusting for $F_2$) is $SS_{1;2}
= Q'_{1;2} (C_{1;2})^- Q_{1;2}$  where $ {\mathbf C}_{1;2}$ and $
{\mathbf Q}_{1;2}$ are as in (\ref  {CQ12}). (See (b) and (c) of
Corollary \ref {rdcedNE}).

Now we seek the answer to the following questions. {\bf Consider a
main effect plan for m factors ($m \geq 3$). Fix a factor, say
$F_i$. What conditions must the design matrices satisfy so that the
sum of squares for $F_i$ adjusted for all others is the same as

(a) the unadjusted sum of squares for $F_i$ ?

(b) the sum of squares for $F_i$ adjusted for only one factor, (say
$F_m$) ?}

[That is so far as $F_i$ is concerned, other factors are virtually
absent.]

\begin{theo} \label{cond-orth} Fix a factor, say $F_i$.

(a)A necessary and sufficient condition for $SS_{i;all} =
SS_{i;m+1}$ is that  the incidence matrix  $N_{ij}$ satisfies the
proportional frequency condition stated in (\ref  {PFC}) [see
Definition \ref {ADDL}.

(b) A necessary and sufficient condition for $SS_{i;all} = SS_{i;m}$
is that
\begin{equation} \label{orthbl-inc} N_{ij} = N_{im} (R_m)^{-1} N'_{jm},
j \neq i,\; 1 \leq i,j \leq m-1. \end{equation}
\end{theo}

The proof relies on two lemmas  we present now.

\begin{lem} \label{P1A&P2A} Consider matrices $A (m \times n), B((m
\times p)$  such that

$$  {\cal C} (B) \subseteq {\cal C} (A) .$$

Let $C((m \times q)$ be any matrix. Then a necessary and sufficient
condition that $ {\cal C} (P_B C) = {\cal C} (P_A C)$ is that $ (P_A
- P_B) C = 0$. \end{lem}

\begin{lem} \label{PA-PD} Consider a matrix $W$ partitioned as
$\left[ \begin{array}{cc} U & V \end{array} \right ]$. Let $ Z = (I
- P_V)U$. Then, $P_W - P_V = P_Z$. \end{lem}

 {\bf Proof of theorem \ref {cond-orth}:} Let $T = \{ 1,2, \cdots i-1,
 i+1, \cdots m\}$ and $T^{*} = T \cup \{m+1 \}$. From
 Notation \ref {sumsofsq} (b), we see that

$$ SS_{i;all} = Y'P_U Y,\; SS_{i;m+1} = Y'P_V Y, $$

 where $U = (I - P_{T^{*}}) X_i$ and  $V = (I - P_{m+1})X_i$.

{\bf Proof of (a);} From the expressions above, a necessary and
sufficient condition for $SS_{i;all} = SS_{i;m+1}$ is that $ P_U =
P_V$, that is ${\cal C} (U) = {\cal C} (V)$.
 Take $ A = X_{m+1}, \: B = X_{T^*}, \: C = X_i$. Then, clearly,
 ${\cal C} (A) \subset {\cal C} (B)$,
that is $ [ {\cal C} (B)]^\bot \subset [{\cal C} (A)]^\bot$. By
Lemma \ref {P1A&P2A} a necessary and sufficient condition for ${\cal
C} (U) = {\cal C} (V)$ is that $[(I - P_A) - (I - P_B)]C =0$, which
is same as

\begin{equation}\label{cond} (P_{T^*} - P_{m+1}) X_i = 0. \end{equation}

Now by Lemma \ref {PA-PD}, $P_{T^*} - P_{m+1} = P_Z$, where $Z = (I
- P_{m+1}) X_T$. Thus, (\ref {cond}) is

$$ \Leftrightarrow P_z X_i = 0\;  \Leftrightarrow X'_i Z = 0
 \;\Leftrightarrow X'_i (I - P_{m+1}) X_j = 0,  j \neq i,$$

which is the same as the proportional frequency condition.

{\bf Proof of (b) :}  Proceeding along similar lines as in the proof
of (a), we find that the necessary and sufficient condition for
$SS_{i;all} = SS_{i;m}$ is that

\begin{equation} P_W X_i = 0,\; \mbox{ where } W = (I - P_m) X_T.
\end{equation}

But this condition $\Leftrightarrow X'_i W = 0 \Leftrightarrow X'_i
(I - P_m) X_j = 0, j \neq i,\; 1 \leq i,j \leq m-1$.  This condition
simplifies to the form in the statement. $\Box$

{\bf Remark 4.4 :} The sufficiency part of (a) of Theorem \ref
{cond-orth} is  well-known. We  now point out that the respective
conditions are also necessary for the sum of squares to satisfy
these desirable properties.

 \vskip5pt

 {\bf Properties of a plan orthogonal through a factor.}

Let us recall Definition \ref {orth3rd}.

\begin{theo} \label{POFm} An MEP orthogonal through $F_m$  has the
following properties.

(a) For every factor $F_i, 1 \leq i \leq m-1$, the reduced normal
equation for $\widehat{\alpha^i}$ is
 $$ ({\mathbf R}_i - {\mathbf N}_{im} (R_m)^{-1} ({\mathbf
 N}_{im})') \; \widehat{\alpha^i} = T_i - {\mathbf N}_{im} (R_m)^{-1} T_m.$$

(b) The error sum of squares is obtained by subtracting the
following from the total sum of squares. Add the sum of squares for
each $F_j$ adjusted for $F_m$, $1 \leq j \leq m-1$ and then the
unadjusted sum of squares for $F_m$. Symbolically,
$$SS_E = SS_{tot} - \sum_{j=1}^{m-1} SS_{j;m} - SS_{m;m+1}.$$
\end{theo}

 {\bf Proof :} Let $L = \{1,2, \cdots m-1 \}$. Then, the
reduced normal equation for the combined effect of the vector of
treatment factors ($\widehat{\alpha^L}$) (after eliminating
$\hat{\mu}$ and $\widehat{\alpha^m}$) is

\begin{eqnarray} \label{blCQall}
 C_{LL;m} \widehat{\alpha^L} & = & Q_{L;m}, \: \mbox{ where  }
 C_{LL;m}  =  ((C_{ij;m}))_{1 \leq i,j \leq m-1}, \ C_{ij;m} =
  X'_i(I - P_{m}) X_j, \\
 Q_{L;m} & = & ((Q_{i;m}))_{1 \leq i \leq m-1}, \; Q_{i;m} =
 X'_i(I - P_m)Y.
                                                                                                                                                                                                                                                                      _{m})Y.\end{eqnarray}

(a) Since the plan is orthogonal through $F_m$,  $C_{ij;m} = 0, i
\neq j, 1 \leq i,j \leq m-1$. Thus,  the reduced normal equation for
$\widehat{\alpha^i}$ is $ C_{ii;m} \widehat{\alpha^i} = Q_{i;m}$,
where $C_{ii;m}$ and $Q_{i;m}$ are as in ( \ref {blCQall}) and the
next equation. That $C_{ii;m}$ and $Q_{i;m}$ are of the form in the
statement of the theorem can be verified easily.

(b) Since the off-diagonal block matrices of $C_{LL;m}$ are null,

$$  SS_{L;m,m+1} =  \sum_{i=1}^{m-1} Q'_{i;m} (C_{ii;m})^- Q_{i;m} =
\sum_{i=1}^{m-1} SS_{i;m}$$

Now the rest follows from (\ref {SSE}). $\Box$

{\bf Remark 4.5:} Statement (a) of Theorem \ref{POFm} was observed
as early as 1996 by Morgan and Uddin in their Theorem 2.1 [equation
(7)]. However, since their paper essentially was concerned with the
construction of optimal nested row-column designs, the result was
overlooked by many authors (such as Mukherjee, Dey and Chatterjee
(2001) and Bagchi (2010)) working on blocked main effect plans.

 {\bf Remark 4.6:} (b) of Theorem \ref{POFm} shows that a plan
orthogonal through one factor ($F_m$) considerably simplifies the
computation of error SS as well as the sum of squares for the
treatment factors $ F_1, F_2 \cdots F_{m-1}$. Thus, in case $F_m$
happens to be a block factor, the whole analysis is only a little
more involved than a fully orthogonal plan, as has been noted in
Bagchi (2010). However, in the situation when $F_m$ is a treatment
factor, analysis is a little more involved since $SS_{m;all}$ needs
to be computed. Needless to mention that the precision of the BLUEs
of the main effect contrasts of $F_m$ (being non-orthogonal to m-1
other factors) is less than the other factors.

  \vskip5pt

{\bf Remark 4.7:} Let us recall the plan $A_5(2)$ [see (\ref
{ORTHthruD})]. If we remove the last column (run) and the last row
(factor D), then we get an OA (4,3,2,2), say $ \rho^*$. Let $C_Q$
denote the coefficient matrix of the reduced normal equation for
factor Q, Q = A,B,C obtained from the plan $ \rho^*$. One may check
that $ C_{Q;\bar{Q}} = C_Q$ for Q = A,B,C. Thus, even though
$A_5(2)$ is not orthogonal, the main effects of factors A,B and C
are estimated with the same precision as the orthogonal plan
$\rho^*$. Therefore, by adding one more run, we are able to
accommodate one more factor (D), without sacrificing the precision
of the three existing factors. The main effect of D is, however,
being estimated with less precision than the others. \vskip5pt

 {\bf  Data analysis of an inter-class orthogonal plan.}

 \vskip5pt

\begin{nota} \label{inter-orth} Consider an inter-class orthogonal plan
with $k$ classes, the ith class having $m_i$ factors denoted by $
F_{i,1}, \cdots F_{i,m_i}, \: 1 \leq i \leq k$. Let $F_G$ denote the
general effect.

 (a) Let $\alpha^{ij}$ and $X^i_j$ denote respectively the vector of
unknown effects  and the design matrix of $F_{ij},\: 1\leq j \leq
m_i,\: i=1, \cdots k$.

 (b) Let $I_i =\{(i,1), \cdots (i,m_i) \}$. For a fixed $j, \:
 1\leq j \leq m_i$, let $T_j = \{(i,j+1), \cdots (i,m_i) \}$ and
 $\bar{j} = I_i \setminus \{(i,j)\}$.

(c) $SS^i_{j;U}$ will denote the sum of squares for $F_{i,j}$
adjusted for each $F_{i,k},\: k \in U$, where $j$ is not in $U$.

(d) $SS^i_{j;all}$ will denote the sum of squares for $F_{i,j}$
adjusted for all other factors in its's own class, i.e.
$SS^i_{j;all} = SS^i_{j;\bar{j}}$.

Further, $SS^i_{j;all>}$ will denote the sum of squares for $F_{ij}$
adjusted for all factors next to it in its's own class, $ 1 \leq j
\leq m_i -1$, while $SS^i_{m_i;all>}$ will denote the  sum of
squares  for $F_{i,m_i}$ adjusted for $F_G$ (that is the unadjusted
sum of squares). Thus,

$$SS^i_{j;all>} = SS^i_{j;T_j}, 1\leq j \leq m_i -1 \mbox{ and } SS^i_{m_i;all>}
= SS_{(i,m_i);G}.$$

 (e) The following expression will be referred to as the
 class total for the ith class.
$$ SS^i_{total} =  \sum_{j=1}^{m_i} SS^i_{j;all>}.$$ \end{nota}

\begin{theo} \label{InclAna} Consider an inter-class orthogonal plan as in Notation
\ref {inter-orth}.  Fix a class, say the ith one and a factor, say
$F_{i,j}$.

 (a) The {\bf reduced normal equation} for $\widehat{\alpha^{ij}}$ is
  {\bf obtained by eliminating only the other factors in the ith class}.
  More explicitly, the reduced normal equation is as follows.
 $$ C^i_{j;\bar{j}} \widehat{ \alpha^{ij}} = Q^i_{j;\bar{j}} \:
  \mbox{ where } C^i_{j;\bar{j}} = (X^i_j)'(I - P^i_{\bar{j}}) X^i_j \:
  \mbox{ and } Q^i_{j;\bar{j}} = (X^i_j)'(I - P^i_{\bar{j}})Y.$$

 Here $P^i_{\bar{j}}$ is the projection operator onto the column
 space of $X^i_{\bar{j}}$.

(b)  The sum of squares for $F_{ij}$, {\bf adjusted for  all other
factors}, is nothing but the sum of squares {\bf adjusted for  all
other factors in the ith class}. Similar statements hold for the sum
of squares adjusted for all factors next to $F_{ij}$. Symbolically,
$$ SS_{(i,j);all} = SS^i_{j;all} \mbox{ and }
SS_{(i,j);all>} = SS^i_{j;all>}.$$

(c) The error sum of squares is obtained by subtracting the class
totals for all the k classes from the total sum of squares.
Symbolically,
 $$SS_E = SS_{tot} -  \sum_{i=1}^{k}  SS^i_{total}.$$
 \end{theo}

   We now note, that if all the factors of a class except one are
   mutually orthogonal through that one, the computation of the class total
   is considerably simpler. The proof follows from  Theorem \ref {POFm}.

\begin{theo} Suppose an inter-class orthogonal plan has a class
(say the ith one) in which all the factors are orthogonal through
$F_{i,m_i}$. Then, the class total for this class can be expressed
as follows.
$$ SS^i_{total} = \sum_{j=1}^{m_i-1} SS^i_{j;m_i} + SS_{m_i;G}.$$
\end{theo}

\section{Main effect plans of small size.}

In this section, we present MEPs with fifteen or less runs obtained
by ad-hoc methods. The factors have at most five levels and the
class size of each plan is at most three. The plans having
class-size three are $A_{12}(1), A_{12}(3)$ and $A_{12}(4)$.
Further, all plans except $A_{12}(2)$ and $A_{12}(3)$ are saturated.
 The graph next to each plan shows the relationship between factors :
the edges drawn with continuous lines represent orthogonality while
dotted line indicate partial orthogonality. The factors are named as
A,B, ....   in the natural order. The equal frequency plans
 are indicated by ``*".

 We begin with a general plan for two $p$-level and one two-level
factors on $2p$ runs. If $p = 3$, the levels of $A$ form a balanced
incomplete block design (BIBD) with those of $B$.

 $A_{2p}(1) = \rho(2p,3; \{p^2\}.2) =$
\begin{eqnarray*}  \left [ \begin{array}{cccccccc}
                 0 & 1 & \cdots & p-1 & 0 & 1 & \cdots & p-1 \\
                 0 & 1 & \cdots & p-1 & 1 & 2 & \cdots & 0 \\
                  0 & 0 & \cdots & 0   & 1 & 1 &\cdots &1 \\
\end{array} \right ] & \setlength{\unitlength}{2mm}
\begin{picture}(20,3)(-2,2)
\thicklines

\put(5,1){\line(-1,1){3}} \put(5,1){\line(1,1){3}}
\put(4.6,0.8){$\bullet$} \put(1.6,3.6){$\bullet$}
\put(7.7,3.6){$\bullet$}

\put(3.3,0.0){C} \put(0.5,3.5){A} \put(8.5,3.5){B}
 \end{picture} \end{eqnarray*}

\vskip5pt

Now a plan with {\bf 6 runs}.

 \begin{center} $ A_6(1) =
\left [\begin{array}{ccccccc}
               A & 0 & 0 & 1 & 1 & 2 & 2 \\
               B & 0 & 1 & 0 & 0 & 1 & 0\\
                C & 0 & 1 & 0 & 1 & 0 & 1 \\
                 D&0 & 1 & 1 & 0 & 0 & 1 \\
\end{array} \right ]$ \end{center}

{\bf Remark 5.1:} For the same experiment an equal-frequency plan is
available - plan $L_6(3.2^3)$ of Wang and Wu (1992). The graphical
representation of these two plans are shown below.
 \begin{eqnarray*} L_6(3.2^3) =
 \setlength{\unitlength}{3mm}
\begin{picture}(6,0)(-1,1)
\thicklines \put(1,1){\line(1,0){3}} \put(1,1){\line(2,1){3}}
\put(1,1){\line(2,-1){3}}

\put(0.7,0.8){$\bullet$} \put(3.8,0.8){$\bullet$}
\put(3.7,2.2){$\bullet$} \put(3.8, -0.8){$\bullet$}

\put(-0.2,0.8){A} \put (4.4,0.6){C} \put(4.4,2.1){B}
\put(4.4,-0.8){D}
\end{picture} \; \; \mbox{ while } A_6(1) =
\setlength{\unitlength}{3mm}
\begin{picture}(20,4)(0,2)

\thicklines

\put(1,1){\line(0,1){3}} \put(1,1){\line(1,1){3}}

 \put(1,4){\line(1,-1){3}} \put(4,1){\line(0,1){3}}
\put(0.7,0.8){$\bullet$} \put(3.8,0.8){$\bullet$}
\put(0.7,3.8){$\bullet$} \put(3.8,3.8){$\bullet$}

\put(-0.2,0.5){C} \put(4.5,0.5){D} \put(-0.2,3.8){A}
\put(4.5,3.8){B} \put(1,3.8){-} \put(2,3.8){-}
\put(3,3.8){-}\put(4,3.8){-}

 \end{picture} \end{eqnarray*}

\vskip5pt

Regarding the performances, the new plan estimates all but one
contrast $(C_1 = \hat{a_0} - 2\hat{a_1} + \hat{a_2}) $ with equal or
more precision. Using the formulae in Theorem \ref {InclAna}, one
may check that the amount of computation is also less here. However,
$C_1$ may be more important for some experimenter, in which case the
old plan $L_6(3.2^3)$ would be preferable.

\vspace{.25em}

We now present plans on  {\bf 8 runs}. We take up the well-known
$OA(8,4,3.2^4)$ and add one more two-level factor (F) with unequal
frequency such that it is orthogonal to all other factors except
$A$.

\vspace{.25em} (a) $A_8(2) = \rho(8,6; \{3\times 2\}.2^4) =$

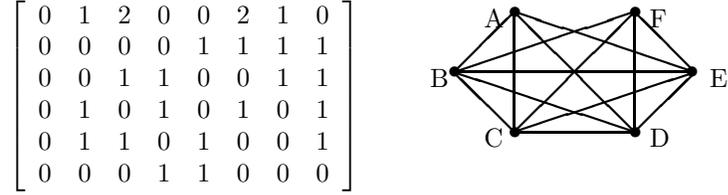
 \begin{eqnarray*}
  \left [\begin{array}{cccccccc}
                             0 & 1 & 2 & 0 & 0 & 2 & 1 & 0 \\
                             0 & 0 & 0 & 0 & 1 & 1 & 1 & 1 \\
                             0 & 0 & 1 & 1 & 0 & 0 & 1 & 1  \\
                             0 & 1 & 0 & 1 & 0 & 1 & 0 & 1 \\
                             0 & 1 & 1 & 0 & 1 & 0 & 0 & 1 \\
                             0 & 0 & 0 & 1 & 1 & 0 & 0 & 0  \\
                              \end{array} \right ] & \; \; \qquad
 \begin{picture}(20,4)(0,34)
\setlength{\unitlength}{4mm} \thicklines \put(2,2){\line(0,1){4}}
\put(2,2){\line(-1,1){2}} \put(2,2){\line(2,0){4}}
\put(2,2){\line(3,1){6}}
 \put(2,2){\line(1,1){4}}

 \put(6,2){\line(-1,1){4}} \put(6,2){\line(-3,1){6}}
 \put(6,2){\line(1,1){2}} \put(6,2){\line(0,1){4}}

 \put(0,4){\line(1,1){2}}\put(0,4){\line(4,0){8}}
 \put(0,4){\line(3,1){6}}

 \put(8,4){\line(-3,1){6}}\put(8,4){\line(-1,1){2}}

\put(1.8,1.8){$\bullet$} \put(5.8,1.8){$\bullet$}
\put(-0.2,3.8){$\bullet$}\put(1.8,5.8){$\bullet$}
\put(7.7,3.8){$\bullet$}\put(5.8,5.8){$\bullet$}

\put(1.0,1.5){C} \put(6.5,1.5){D} \put(-0.8,3.5){B} \put(1.0,5.5){A}
\put(8.5,3.5){E} \put(6.5,5.5){F}

 \end{picture}
 \end{eqnarray*}

Our next {\bf plan on 8 runs} has two three-level factors satisfying
partial orthogonality.

(b) The plan  $A_8(3) = \rho(8,5; \{3^2\}.\{2^2\}.2) =$
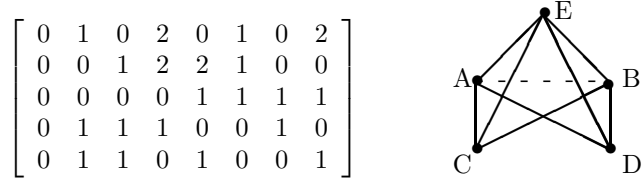
\begin{eqnarray*}
 \left [ \begin{array}{cccccccc}
           0 & 1 & 0 & 2 & 0 & 1 & 0 & 2 \\
          0 & 0 & 1 & 2 & 2 & 1 & 0 & 0 \\
          0 & 0 & 0 & 0 & 1 & 1 & 1 & 1 \\
          0 & 1 & 1 & 1 & 0 & 0 & 1 & 0 \\
          0 & 1 & 1 & 0 & 1 & 0 & 0 & 1 \\
           \end{array} \right ] & \; \;\qquad
           \setlength{\unitlength}{3mm}
\begin{picture}(20,6)(0,3) \thicklines
\put(1,1){\line(0,1){3}} \put(1,1){\line(2,1){6}}
 \put(1,1){\line(1,2){3}}
 \put(1,4){\line(2,-1){6}} \put(1,4){\line(1,1){3}}
  \put(4,7){\line(1,-2){3}}  \put(4,7){\line(1,-1){3}}
 \put(7,1){\line(0,1){3}} \put(7,1){\line(-1,2){3}}

 \put(0.8,0.8){$\bullet$} \put(0.8,3.8){$\bullet$} \put(3.8,6.8){$\bullet$}
 \put(6.7,0.8) {$\bullet$} \put(6.6,3.6){$\bullet$}

\put(0.0,0.0){C} \put(0.0,3.8){A} \put(4.5,6.8){E} \put(7.5,0.0){D}
\put(7.5,3.8){B} \put(1,3.8){-} \put(2,3.8){-}
\put(3,3.8){-}\put(4,3.8){-}
 \put(5,3.8){-}\put(6,3.8){-}
\end{picture}  \end{eqnarray*}

\vspace{.5em}

 We now present two plans with 4-level factors on 8
runs. Note that on 8 runs, a four-level factor can be orthogonal to
neither a four-level nor a three-level factor. Using
non-orthogonality, we are able to accommodate two four-level factors
in one plan and one four-level and one three-level factor in another
plan on 8 runs. \vspace{.5em}

(c)  $A_8(4)^* = \rho(8,3; \{4^2\}.2) =$ is obtained by putting
$p=4$ in $A_{2p}(1)$. In this plan,  the 4-level factors A and B
form a group divisible design $(m=n=2, r=k=2, \lambda_1 =0,
\lambda_2=1)$.

\vspace{.5em}

(d)  $A_8(5) = \rho(8,4; \{4.3\}.2^2) =$
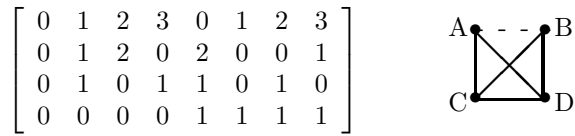
\begin{eqnarray*}   \left [ \begin{array}{cccccccc}
                  0 & 1 & 2 & 3 &  0 & 1 & 2 & 3  \\
                  0 & 1 & 2 & 0 &  2 & 0 & 0 & 1 \\
                  0 & 1 & 0 & 1 &  1 & 0 & 1 & 0 \\
                  0 & 0 & 0 & 0 & 1 & 1 & 1 & 1 \\
                  \end{array} \right ]& \; \; \qquad
 \setlength{\unitlength}{3mm}
\begin{picture}(20,4)(0,2)

\thicklines

\put(1,1){\line(1,0){3}} \put(1,1){\line(0,1){3}}
\put(1,1){\line(1,1){3}}

 \put(1,4){\line(1,-1){3}} \put(4,1){\line(0,1){3}}
\put(0.7,0.8){$\bullet$} \put(3.8,0.8){$\bullet$}
\put(0.7,3.8){$\bullet$} \put(3.8,3.8){$\bullet$}

\put(-0.2,0.5){C} \put(4.5,0.5){D} \put(-0.2,3.8){A}
\put(4.5,3.8){B} \put(1,3.8){-} \put(2,3.8){-}
\put(3,3.8){-}\put(4,3.8){-}

 \end{picture} \end{eqnarray*}
\vspace{.5em}

{\bf A plan on 10 runs :}  $A_{10}^{*} = \rho(10,3; \{5^2\}.2) =$ is
obtained by putting $p=5$ in $A_{2p}(1)$. Here the 5-level factors
form a symmetric cyclic PBIBD with $r=k=2, \lambda_1 =1,
\lambda_2=0)$.

\vspace{.5em}

We shall now present {\bf plans on 12 runs :}. There is no plan in
the literature accommodating one or more 4-level factors on 12 runs.
So, we begin with such plans.

(a)  $A_{12}(1)^* =  \rho(12,5; \{2.4^2 \}.\{3^2\}) =$
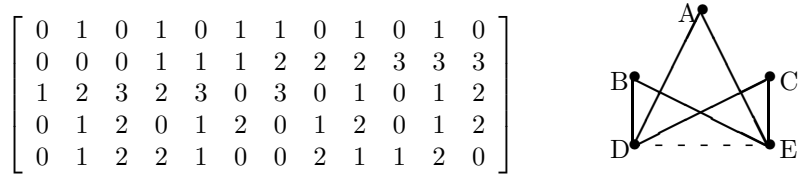
\begin{eqnarray*}  \left [ \begin{array}{cccccccccccc}

        0 & 1 & 0 & 1 & 0 & 1 & 1 & 0 & 1 & 0 & 1 & 0 \\
        0 & 0 & 0 & 1 & 1 & 1 & 2 & 2 & 2 & 3 & 3 & 3 \\
        1 & 2 & 3 & 2 & 3 & 0 & 3 & 0 & 1 & 0 & 1 & 2 \\
        0 & 1 & 2 & 0 & 1 & 2 & 0 & 1 & 2 & 0 & 1 & 2   \\
        0 & 1 & 2 & 2 & 1 & 0 & 0 & 2 & 1 & 1 & 2 & 0 \\
                            \end{array} \right ] & \; \; \qquad
 \setlength{\unitlength}{3mm}
\begin{picture}(20,4)(0,3)

\thicklines

 \put(1,1){\line(0,1){3}} \put(1,1){\line(2,1){6}}
 \put(1,1){\line(1,2){3}}

 \put(1,4){\line(2,-1){6}} \put(4,7){\line(1,-2){3}}
 \put(7,1){\line(0,1){3}}

 \put(0.8,0.8){$\bullet$} \put(0.8,3.8){$\bullet$} \put(3.8,6.8){$\bullet$}
 \put(6.8,0.8) {$\bullet$} \put(6.8,3.8){$\bullet$}

 \put(0.0,0.5){D}  \put(0.0,3.5){B}  \put(3.0,6.5){A}
 \put(7.5,0.5){E}  \put(7.5,3.5){C}

\put(1,0.8){-} \put(2,0.8){-}\put(3,0.8){-}\put(4,0.8){-}
\put(5,0.8){-} \put(6,0.8){-}

 \end{picture} \end{eqnarray*}
 \vspace{.5em}

In this equal frequency saturated plan, both the four-level factors
B and C  form a generalized group divisible design with the levels
of the two-level factor A. Between themselves, they form a Balanced
incomplete block design (BIBD). The relation between factors $D$ and
$E$ is presented in details after Remark 2.4.

\vspace{.5em}

(b) $A_{12}(2) = \rho(12,4; \{3^2\}.\{3.4\}) =$
\begin{eqnarray*}  \left [
\begin{array}{cccccccccccc}
         0 & 0 & 0 & 0 & 1 & 1 & 1 & 1 & 2 & 2 & 2 & 2 \\
         0 & 1 & 2 & 0 & 1 & 2 & 0 & 1 & 2 & 0 & 1 & 2   \\
         0 & 1 & 2 & 0 & 0 & 0 & 1 & 2 & 1 & 2 & 0 & 0 \\
         0 & 1 & 2 & 3 & 0 & 1 & 2 & 3 & 0 & 1 & 2 & 3 \\
  \end{array} \right ] & \; \; \qquad
 \setlength{\unitlength}{4mm}
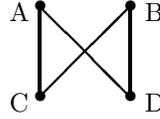
\begin{picture}(20,3)(0,2)

\thicklines

 \put(1,1){\line(0,1){3}} \put(1,1){\line(1,1){3}}

 \put(1,4){\line(1,-1){3}} \put(4,1){\line(0,1){3}}
\put(0.8,0.8){$\bullet$} \put(3.8,0.8){$\bullet$}
\put(0.8,3.8){$\bullet$} \put(3.8,3.8){$\bullet$}

\put(0.0,0.5){C} \put(4.5,0.5){D} \put(0.0,3.5){A} \put(4.5,3.5){B}

 \end{picture} \end{eqnarray*}

Here all the factors other except C has equal frequency. The levels
of factors A and B form a balanced block design (BBD). $D$ is
partially orthogonal to $C$ as contrast $\hat{\delta_1} -
\hat{\delta_2}$ is orthogonal the contrasts for $C$. However, $C$ is
non-orthogonal to $D$.

\vspace{.5em}

 {\bf Remark 5.2:} In the plan $A_{12}(2)$, the
four-level factor D may be replaced by three mutually orthogonal
two-level factors to obtain an almost orthogonal MEP for an $3^3.
2^3$ experiment. \vspace{.5em}

 (c) $A_{12}(3)^* = \rho(12,7; 2^4.\{3^3\}) =$
\begin{eqnarray*} \left [ \begin{array}{cccccccccccc}
                    0 & 0 & 1 & 1 & 0 & 0 & 1 & 1 & 0 & 0 & 1 & 1  \\
                    0 & 1 & 0 & 1 & 0 & 1 & 0 & 1 & 0 & 1 & 0 & 1 \\
                    0 & 0 & 1 & 1 & 1 & 1 & 0 & 0 & 0 & 1 & 1 & 0 \\
                    0 & 1 & 1 & 0 & 0 & 1 & 0 & 1 & 1 & 0 & 1 & 0 \\
                    0 & 0 & 0 & 0 & 1 & 1 & 1 & 1 & 2 & 2 & 2 & 2  \\
                    0 & 1 & 0 & 2 & 1 & 2 & 2 & 0 & 2 & 0 & 1 & 1 \\
                    0 & 1 & 1 & 2 & 2 & 0 & 1 & 2 & 2 & 1 & 0 & 0\\
   \end{array} \right ] & \; \; \qquad
 \setlength{\unitlength}{4mm}
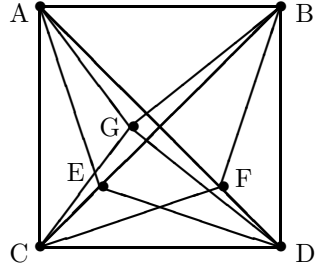
\begin{picture}(20,3)(0,4)
\thicklines

 \put(1,1){\line(8,0){8}} \put(1,1){\line(0,8){8}}
\put(1,1){\line(1,1){8}} \put(1,1){\line(1,1){2}}
\put(1,1){\line(3,1){6}} \put(1,1){\line(3,4){3}}

\put(9,1){\line(-1,1){8}} \put(9,1){\line(0,8){8}}
\put(9,1){\line(-3,1){6}} \put(9,1){\line(-1,1){2}}
\put(9,1){\line(-5,4){5}}

\put(1,9){\line(8,0){8}} \put(1,9){\line(1,-3){2}}
 \put(1,9){\line(1,-1){6}}
\put(1,9){\line(3,-4){3}}

\put(3,3){\line(1,1){6}} \put(7,3){\line(1,3){2}}
 \put(4,5){\line(5,4){5}}

\put(0.8,0.8){$\bullet$} \put(8.8,0.8){$\bullet$}
\put(0.8,8.8){$\bullet$} \put(8.8,8.8){$\bullet$}
\put(2.9,2.8){$\bullet$} \put(6.9,2.8){$\bullet$}
\put(3.9,4.8){$\bullet$}

\put(0.0,0.5){C} \put(9.5,0.5){D} \put(0.0,8.5){A} \put(9.5,8.5){B}
\put(1.9,3.2){E} \put(7.5,3.0){F} \put(3.0,4.7){G}

\end{picture} \end{eqnarray*}

\vspace{.5em}

{\bf Remark 5.3:} The plan  $A_{12}(3)$ is very similar to the plan
$L_{12}^{\prime}(3^4.2^3)$ of Wang and Wu (1992). The difference is
that $A_{12}(3)$ provides  one more two-level factor and one less
three-level factor and so has total d.f one less than
$L_{12}^{\prime}(3^4.2^3)$. On the other hand, since  each
three-level factor (say P) in $A_{12}(3)$ is non-orthogonal to two
and not three factors and the relationship of P with the any other
three-level factor is the same as that in $L_{12}^{\prime}(3^4.2^3)$
it's contrasts are estimated with greater precision.

\vspace{.5em}

(d) $A_{12}(4)= \rho(12,6; \{3^3\}.\{3^2.2\}) =$
\begin{eqnarray*} \left [ \begin{array}{cccccccccccc}
                    0 & 0 & 0 & 0 & 1 & 1 & 1 & 1 & 2 & 2 & 2 & 2  \\
                    2 & 1 & 2 & 0 & 0 & 2 & 1 & 1 & 0 & 1 & 0 & 2 \\
                    0 & 1 & 2 & 1 & 0 & 2 & 1 & 2 & 1 & 2 & 0 & 0\\
                    0 & 0 & 1 & 2 & 0 & 0 & 1 & 2 & 0 & 0 & 1 & 2 \\
                     2 & 1 & 0 & 0 & 0 & 1 & 0 & 2 & 2 & 0 & 1 & 0\\
                      0 & 1 & 1 & 0 & 1 & 0 & 0 & 1 & 1 & 0 & 0 & 1\\
   \end{array} \right ]& \; \; \qquad
\setlength{\unitlength}{4mm}
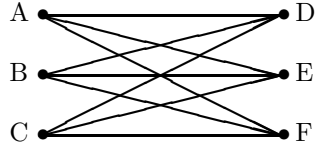
\begin{picture}(20,3)(0,2)
\thicklines

\put(1,4){\line(2,-1){8}} \put(1,4){\line(4,-1){8}}
\put(1,4){\line(1,0){8}}

\put(1,2){\line(4,-1){8}}\put(1,2){\line(4,1){8}}
\put(1,2){\line(1,0){8}}

\put(1,0){\line(4,1){8}}\put(1,0){\line(2,1){8}}
\put(1,0){\line(1,0){8}}

\put(0.9,-0.2){$\bullet$}\put(8.9,-0.2){$\bullet$}
\put(0.9,1.8){$\bullet$} \put(8.9,1.8){$\bullet$}
\put(0.9,3.8){$\bullet$} \put(8.9,3.8){$\bullet$}

\put(0.0,-0.2){C} \put(9.5,-0.2){F}

\put(0.0,1.8){B}\put(9.5,1.8){E}

\put(0.0,3.8){A}\put(9.5,3.8){D}


\end{picture} \end{eqnarray*}

\vspace{.5em}

This inter-class (3) orthogonal MEP  has accommodated five
three-level factors together with a two-level factor. The levels of
A form a BBD with those of each of B and C, while the levels of $B$
form a variance-balanced non-binary design with those of $C$. Both
the three-level factors $D$ and $E$ are partially orthogonal to the
two-level factor $F$.\vspace{.5em}

{\bf A plan on 15 runs :}
 $A_{15} =  \rho(15,4; \{3^2\}.\{5^2\})=$
 \begin{eqnarray*} \left [ \begin{array}{ccccccccccccccc}
    0 & 1 & 2 & 0 & 1 & 2 & 0 & 1 & 2 & 0 & 1 & 2 & 0 & 1 & 2\\
       0 & 1 & 2 & 2 & 0 & 1 & 1 & 2 & 0 & 0 & 1 & 2 & 2 & 0 & 1\\
       0 & 0 & 0 & 1 & 1 & 1 & 2 & 2 & 2 & 3 & 3 & 3 & 4 & 4 & 4\\
       0 & 1 & 2 & 1 & 2 & 3 & 2 & 3 & 4 & 3 & 4 & 0 & 4 & 0 & 1 \\
      \end{array} \right ] & \; \; \qquad
 \setlength{\unitlength}{4mm}
\begin{picture}(20,3)(0,2)

\thicklines

 \put(1,1){\line(0,1){3}} \put(1,1){\line(1,1){3}}

 \put(1,4){\line(1,-1){3}} \put(4,1){\line(0,1){3}}

\put(0.8,0.8){$\bullet$} \put(3.8,0.8){$\bullet$}
\put(0.8,3.8){$\bullet$} \put(3.8,3.8){$\bullet$}

\put(0.0,0.5){C} \put(0.0,3.5){A} \put(4.5,0.5){D} \put(4.5,3.5){B}
 \end{picture} \end{eqnarray*}
\vskip 10pt

\section{References}
\begin{enumerate}
\item Addleman,S. (1962). Orthogonal main effect plans for asymmetrical
factorial experiments.Technometrics 4, p: 21-46.

\item  Bagchi, S. (2010). Main effect plans orthogonal through the
block factor, technometrics, vol. 52, p : 243-249.

\item Dey, A. (1985). Orthogonal fractional factorial designs,
John wiley, New York.

\item  Dey, A. and Mukherjee, R. (1999). Fractional factorial plans,
Wiley Series in Prob. and Stat.

\item Hedayat, A.S., Sloan, N.J.A. and Stufken,J. (1999). Orthogonal
arrays, Theory and Applications, Springer Series in Statistics.

\item  Huang, L., Wu, C.F.J. and Yen, C.H. (2002). The idle
column method : Design construction, properties and comparisons,
Technometrics, vol. 44, p : 347-368.

\item  Ma, C. X., Fang, K.T. and Liski, E. (2000). A new approach in
constructing orthogonal and nearly orthogonal arrays, Metrika, vol.
50, p : 255-268.

\item Morgan, J.P. and Uddin, N. (1996). Optimal blocked main efect
plans with nested rows and columns and related designs, Ann. Stat.,
vol. 24, p : 1185-1208.

\item  Mukherjee, R., Dey, A. and Chatterjee, K. (2001). Optimal
main effect plans with non-orthogonal blocking. Biometrika, vol. 89,
p : 225-229.

\item  Nguyen, (1996). A note on the construction of
near-orthogonal arrays with mixed levels and economic run size,
Technometrics, vol. 38, p : 279-283.

\item Starks, T.H. (1964). A note on small orthogonal main
effect plans for factorial experiments, Technometrics, 8, P :
220-222.

\item Wang, J.C. and Wu.,C.F.J. (1992). Nearly orthogonal arrays with
mixed levels and small runs, Technometrics, vol. 34, p : 409-422.

\item Xu, H. (2002). An algorithm for constructing orthogonal and nearly orthogonal
arrays with mixed levels and small runs, Technometrics, vol. 44, p :
356-368.
\end{enumerate}

\end{document}